\documentclass{article}
\usepackage{amssymb}
\usepackage{amsmath}
\usepackage{amsthm}

\usepackage{caption}

\usepackage[colorlinks,citecolor=blue,urlcolor=blue,linkcolor=blue,pagecolor=blue]{hyperref}

\usepackage{enumerate}
\usepackage{pgf,pgfarrows,pgfnodes}

\newcommand{\eq}{\begin{equation}}
\newcommand{\en}{\end{equation}}
\newcommand{\re}[1]{\mbox{(\ref{#1})}}

\newtheorem{Theorem}{Theorem}
\newtheorem{theorem}[Theorem]{Theorem}
\newtheorem{lemma}[Theorem]{Lemma}
\newtheorem{corollary}[Theorem]{Corollary}
\newtheorem{construction}[Theorem]{Construction}
\newtheorem{proposition}[Theorem]{Proposition}
\newtheorem{example}[Theorem]{Example}
\newtheorem{defn}[Theorem]{Definition}
\newtheorem*{defn*}{Definition}

\newtheorem{question}[Theorem]{Question}
\newtheorem{conjecture}[Theorem]{Conjecture}
\newtheorem{condition}[Theorem]{Condition}
\newtheorem{remark}[Theorem]{Remark}
\newtheorem{problem}[Theorem]{Problem}
\def\proof{\noindent{\bf Proof.\ \ }}

\def\endpf{\hfill $\square$ \vskip .25in}

\renewcommand{\P}{\mathbb{P}}

\newcommand{\E}{\mathbb{E}}

\newcommand{\te}{\rightarrow}
\newcommand{\ed}{\mbox{$ \ \stackrel{d}{=}$ }}

\newcommand{\lb}[1]{\label{#1}}

\newcommand {\giv}{ \,|\,}

\def\proof{\noindent{\bf Proof.\ \ }}

\def\endpf{$\Box$}

\begin{document}

\title{Concave Majorants of Random Walks and Related Poisson Processes}
\author{Josh Abramson\thanks{University of California at Berkeley;  josh@stat.berkeley.edu}
        \ and 
        Jim Pitman\thanks{University of California at Berkeley;  pitman@stat.berkeley.edu; research supported in part by the National Science Foundation Award $0806118$} 
      }
\date{\today}
\maketitle 

\begin{abstract}
We offer a unified approach to the theory of concave majorants of random walks 
by providing a path transformation for a walk of finite length that leaves the law of the walk unchanged 
whilst providing complete information about the concave majorant. 
This leads to a description of a walk of random geometric length as a Poisson point process of excursions away from its concave majorant, 
which is then  used to find a complete description of the concave majorant of a walk of infinite length. 
In the case where subsets of increments may have the same arithmetic mean, we investigate three nested compositions that naturally arise from our construction of the concave majorant.
\end{abstract}

\section{Introduction}
\label{sec:intro}

Let $S_0 = 0$ and $S_j = \sum_{i=1}^j X_i$ for $1 \le j \le n$, where $X_1, \ldots , X_n$ 
are exchangeable random variables. Let $\mathbf{A}$ be the assumption that almost surely no two subsets of $X_1,\ldots, X_n$ have the same arithmetic mean, and assume for now that $\mathbf{A}$ holds. 
Let $S^{[0,n]} := \{ (j,S_j) \: : \: 0 \le j \le n \}$, so that $S^{[0,n]}$ is the random walk of length $n$ with increments distributed like $X_1,\ldots,X_n$. 
Let
$$
0 < N_{n,1} < N_{n,1} + N_{n,2} < \cdots < N_{n,1} + \cdots + N_{n,F_n} = n
$$
be the successive times $j$ with $0 \le j \le n$ such that $S_j = \bar{C}^{[0,n]}(j)$,
where $ \bar{C}^{[0,n]}$
is the {\em concave majorant} of the walk $S^{[0,n]}$, i.e.\ the least
concave function $C$ on $[0,n]$ such $C(j) \ge S_j$ for $1 \le j \le n$. The random variable
$F_n$ is the {\em number of faces} of the concave majorant. Without assumption $\mathbf{A}$, more care needs to be taken in defining the faces of the concave majorant; this will be discussed further in Section \ref{sec:noncontinuous}.

The $i$th face of the concave majorant is a chord from $(N_{n,1}+ \cdots + N_{n,i-1},S_{N_{n,1}+ \cdots + N_{n,i-1} })$ to $(N_{n,1}+ \cdots + N_{n,i},S_{N_{n,1}+ \cdots + N_{n,i} })$. We define the {\em length}, {\em increment} and  {\em slope} of the $i$th face to be $N_i$, $\Delta_{n,i}$ and $\frac{\Delta_{n,i}}{N_i}$ respectively, where
$$
\Delta_{n,i} := ( S_{N_{n,1} + \cdots + N_{n,i}} -  S_{N_{n,1} + \cdots + N_{n,i-1}} ), \quad \text{ for }1 \le i \le F_n.
$$

In the 1950's, E. Sparre Andersen \cite{MR0058893b} discovered the following
remarkable result: for any exchangeable $X_1,\ldots,X_n$ satisfying assumption $\mathbf{A}$,
there is the equality in distribution
\eq
\lb{fnkn}
F_n \ed K_n = \sum_{j = 1}^n I_j
\en
where $K_n$ is the number of cycles in a uniformly distributed random permutation
of the set $[n]:= \{1, \ldots , n\}$, and $I_j, j = 1, 2,\ldots$ is a sequence of independent
Bernoulli variables with $\P(I_j = 1) = 1/j$ and $\P(I_j = 0) = 1- 1/j$  for each $j$.
The second equality in  \re{fnkn} is an elementary and well known representation of $K_n$ which 
holds for a number of natural constructions of uniform random permutations of $n$  simultaneously 
for all $n$, including both the construction from records of the $X_i$ \cite{MR994088},  and the
Chinese Restaurant Process \cite{MR2245368}.

A further result that seems to have been known by Spitzer \cite{MR0079851}, and shown explicitly by Goldie \cite{MR994088}  using a generalization by Brunk of Spitzer's Lemma \cite{MR0162302}, is that  under assumption $\mathbf{A}$ the distribution of 
the {\em partition of $n$} generated by the lengths of the faces of the concave majorant on $[0,n]$,
which may be encoded by these lengths in non-increasing order, 
has the same distribution as the partition of $n$ generated by the cycles of  a uniform random
permutation - we will prove this result as a corollary of our main theorem.
Thus the partition generated by the lengths of the faces of the concave majorant may be generated by a discrete \emph{uniform stick breaking process} on $[0,n]$ \cite{MR2245368}. The result raises the following problem:

\paragraph{The rearrangement problem.} Conditionally given that the partition of $n$ generated by
the lengths of the faces of the concave majorant of the random walk $S^{[0,n]}$ has segment lengths $n_1, \ldots, n_k$ with
$n_1 \ge n_2 \ge \ldots \ge n_k > 0$,
\begin{itemize}
\item
in what order and with what increments should the faces $f_1, \ldots , f_k$ of the concave majorant with lengths $n_1, \ldots , n_k$ respectively be arranged to recreate the concave majorant of the random walk $S^{[0,n]}$?
\item
given the concave majorant, what is the distribution of values of the random walk $S^{[0,n]}$ between vertices of the concave majorant?
\end{itemize}
We answer this question by giving in Theorem \ref{theorem:1}  a simultaneous construction of the walk and its concave majorant conditional on the partition generated by the lengths of the faces of the concave majorant. The theorem will be proved  under assumption $\mathbf{A}$ in Section \ref{sec:proof}, and in the general case in Section \ref{sec:noncontinuous}, with the key idea of both proofs being that it is enough to show that the theorem is true when $X_1,\ldots,X_n$ are samples without replacement from a set of $n$ real numbers. Since the construction given in the theorem applies to general exchangeable $X_1,\ldots,X_n$ it allows us to investigate in Section \ref{sec:noncontinuous} the structure of the concave majorant in the general case.
The statement of the theorem is complicated, but easy to describe informally, particularly under assumption $\mathbf{A}$, in which case the construction is as follows. Conditional on the lengths of the blocks of the partition generated by the concave majorant being $(n_1,\ldots,n_k)$:
\begin{itemize}
\item
Split $X_1,\ldots,X_n$ into $k$ blocks
$$
(X_1,\ldots,X_{n_1})(X_{n_1+1},\ldots,X_{n_1+n_1}) \cdots (X_{\sum_{i=1}^{k-1}n_i + 1} ,\ldots,X_{\sum_{i=1}^k n_i })
$$
\item
Arrange the blocks in order of decreasing arithmetic means.
\item
Perform the unique cyclic permutations of the increments within each block such that the walk with those cyclically permuted increments remains below the line joining its start and end points.
\end{itemize}
This process defines a permutation of the original increments which leaves the distribution of the walk $S^{[0,n]}$ unchanged and at the same time provides us with information about the concave majorant. In the case where $X_1,\ldots,X_n$ are independent, then we may just generate independent walks of length $n_1, \ldots, n_k$, cyclically permute the increments of each walk appropriately, and then arrange the walks in order of decreasing slope. The idea of using cyclic permutations to transform random walk bridges into excursions is due to Vervaat \cite{MR515820}.

When assumption $\mathbf{A}$ is not satisfied there are two more complications. Some of the blocks may have the same arithmetic mean, in which case their ordering is chosen uniformly, and within a block there may be more than one cyclic permutation of increments that leaves the walk with those increments below the line joining its start and end points, in which case the cyclic permutation is chosen uniformly from the possible options. By exchangeability, it would also work to take the blocks with the same arithmetic mean in order of appearance rather than randomly ordering them, but this makes the statement of the theorem harder and in fact does not make the proof any easier.

To facilitate the statement of the theorem, it is necessary to define the set of all permutations that cyclically permute increments within certain blocks and then arrange those blocks in some order.

\begin{defn*}
Let $\Sigma_n$ be the set of permutations of $[n]$, and let $\mathcal{P}_n$ be the set of partitions of $n$, encoded in non-increasing order. For $(n_1,\ldots,n_k) \in \mathcal{P}_n$ let $\Sigma_{(n_1,\ldots,n_k)} \subseteq \Sigma_n$ be such that $\sigma \in \Sigma_{(n_1,\ldots,n_k)}$ if and only if for some $\tau \in \Sigma_k$ and $(r_1,\ldots,r_k) \in \mathbb{Z}^k$ we have
$$
\sigma \left( \mbox{$\sum_{l=1}^{i-1} n_{\tau(l)}$} \: + \: j \right)
= \mbox{$\sum_{l=1}^{\tau(i)-1} n_l$} + \left( (j-1+r_i) \! \! \! \! \mod n_{\tau(i)} \right) + 1
$$
for $1 \le j \le n_{\tau(i)}$, $1 \le i \le k$.
\end{defn*}

In the definition of $\Sigma_{(n_1,\ldots,n_k)}$ just given, the cyclic shift chosen for the $\tau(i)$th block is given by $r_i$ and the ordering of the $k$ blocks is given by $\tau$.

\begin{theorem}
\label{theorem:1}
Let $S_0 = 0$ and $S_j = \sum_{\ell =1}^j X_{\ell}$ for $1 \le j \le n$, where $X_1, \ldots , X_n$ 
are random variables with any exchangeable joint distribution. Let $S^{[0,n]} = \{ (j,S_j) \: : \: 0 \le j \le n \}$.
Independently of $X_1,\ldots,X_n$, let $L_{n,1}, L_{n,2} , \ldots ,  L_{n,K_n}$ be a sequence of random variables distributed like the lengths of cycles of 
a random permutation of $[n]$ arranged in non-increasing order.
Conditionally given $\{ K_n = k \}$ and $\{ L_{n,i} = n_i \: : \: 1
\le i \le k \}$, let $B$ be the random subset of $\Sigma_n$ defined by the following relation. $\sigma$ is in $B$ if and only if $\sigma \in \Sigma_{(n_1,\ldots,n_k)}$ and there exists $\tau \in \Sigma_k$ such that the function defined on $[k]$ by
\eq
\lb{decreasing}
i \mapsto \Delta_{n,i}^{\sigma,\tau} := \frac{1}{ n_{\tau(i)} } \left( \sum_{ \ell =  n_{\tau(1)} + \cdots + n_{\tau(i-1)} + 1 }^{  n_{\tau(1)} + \cdots + n_{\tau(i)} } X_{\sigma(\ell)} \right) 
\en
is non-increasing in $i$ and for each $1 \le i \le k$ we have
\eq
\lb{below}
\frac{1}{ m } \left( \sum_{ \ell =  n_{\tau(1)} + \cdots + n_{\tau(i-1)} + 1 }^{  n_{\tau(1)} + \cdots + n_{\tau(i-1)} + m} X_{\sigma(\ell)} \right) \le \Delta_{n,i}^{\sigma,\tau} \quad \text{ for } 1 \le m \le n_{\tau(i)} . 
\en
Conditionally given $B$, let $\rho$ be a uniform random element of $B$, independently of all previously introduced random variables. For $1 \le j \le n$ let $S_j^{\rho} = \sum_{\ell=1}^j X_{\rho( \ell )}$ and let $S_{\rho}^{[0,n]} = \{ (j,S_j^{\rho}) \: : \: 0 \le j \le n \}$.
Then $S_{\rho}^{[0,n]} \ed S^{[0,n]}$.
\end{theorem}

The condition involving \re{decreasing} ensures that the permutation
that we end up choosing puts the blocks of increments in
non-increasing order of arithmetic mean, i.e.\ in non-increasing order of slope, and the condition involving \re{below} ensures that the cyclic permutation chosen for each block makes the walk stay below the line joining the start and end points of the increments of that block. In the case where $X_1,\ldots,X_n$ satisfy assumption $\mathbf{A}$, the random set $B$ almost surely only consists of one element and thus the additional random variable $\rho$ is not needed.


Some of the ideas of our construction are contained within the work of Spitzer \cite{MR0079851}, who observed that if 
$\Delta_{n,i} $
is the increment of the walk over the $i$th face of the concave majorant, then for the
maximum
$$
M_n := \max_{0 \le k \le n} S_k  
$$
there is the almost sure representation
\eq
\lb{spitzer1}
M_n = \sum_{i = 1}^{F_n} \Delta_{n,i} 1 ( \Delta_{n,i} \ge 0 ).
\en
Spitzer showed the much simpler representation in distribution
\eq
\lb{spitzer2}
M_n \ed \sum_{i = 1}^{K_n} \Delta_{n,i} ^* 1 ( \Delta_{n,i} ^* \ge 0 )
\en
where $K_n$ is the number of cycles of a random permutation independent of the random walk 
$S^{[0,n]} = \{ (j,S_j) \: : \: 0 \le j \le n \} $, and given $K_n = k$ and that the permutation has cycles of lengths 
say $L_{n,1}, \ldots, L_{n,k}$,
the $\Delta_{n,i}^*$ are conditionally independent, with 
$$
( \Delta_{n,i}^*  \giv K_n = k, L_{n,i} = \ell )  \ed S_{\ell} , \quad \text{ for } 1 \le i \le k , \text{ and } 1 \le \ell \le n.
$$
This is an immediate corollary of our theorem, and something we investigate further in Section \ref{sec:max}.
Some consequences of this result lead to other ideas which arise in this paper.
Let $S_\ell^+ = S_\ell \vee 0$. As pointed out by Spitzer, Hunt's remarkable identity \cite[Theorem $4.1$]{MR0062867}
\eq
\lb{hunt}
\E( M_n ) = \sum_{\ell = 1}^n \frac{ \E(S_\ell^+) }{\ell}
\en
follows easily from \re{spitzer2}, along with the following complete 
description of the distribution of $M_n$ for every $n = 1,2, \ldots$ (this description is known as Spitzer's Identity):
for $|q| <  1$
\eq
\lb{spitzer}
\sum_{n = 0}^\infty q^{n} \E e^{i t M_n } = \exp \left( \sum_{k = 1}^\infty \frac{q^k}{k} \E e^{it S_k^+}\right)
\en
To indicate how \re{hunt} follows from \re{spitzer2}, recall 
that the expected number of cycles of length $\ell$ in a random permutation of $[n]$ is $\ell^{-1}$.
So \re{hunt} decomposes the expectation of the sum in \re{spitzer2} according the contributions
from cycles of various sizes $\ell$. To provide a similar interpretation of \re{spitzer},
let $n(q)$ denote a random variable with geometric distribution with parameter $1-q$, so
$\P(n(q) \ge n) = q^n$ for $n = 0,1, \ldots$, and assume $n(q)$ is independent of the random walk. 
Then multiplying \re{spitzer} by $1-q$ and using the
expansion $-\log(1-q) = \sum_{k=1}^\infty q^k/k$ allows \re{spitzer} to be rewritten \cite{greenpit80}:
\eq
\lb{spitzer3}
\E e^{i t M_{n(q)} } = \exp \left( \sum_{k = 1}^\infty \frac{q^k}{k} ( \E e^{it S_k^+ } - 1) \right)
\en
Otherwise put, the maximum $M_{n(q)}$ of the walk up to the independent geometric time $n(q)$
has a compound Poisson distribution:
\eq
\lb{mnq}
M_{n(q)} \ed \sum_{k=1}^\infty \sum_{i=1}^{N(q^k/k)} S_{k,i}^+
\en
where for fixed $q$ the $N(q^k/k)$ are independent Poisson variables with parameters $q^k/k$
for $k = 1,2, \ldots$, and given these variables the $S_{k,i}$ for $1 \le i \le N(q^k/k)$ are
independent with $S_{k,i} \ed S_k$. As observed by Greenwood and Pitman \cite{greenpit80}, the identity
in distribution \re{mnq}, and the companion result which determines the common distribution of 
$S_{n} - M_{n}$ and $\min_{0 \le k \le n}S_k$ for every $n$, can be derived, along with other
results of fluctuation theory for the distribution of ladder heights and ladder times, from the
decomposition 
\eq
S_{n(q)} = M_{n(q)} + ( S_{n(q)} - M_{n(q)} )
\en
which expresses the compound Poisson variable $S_{n(q)}$ as the sum of two independent 
compound Poisson variables with with positive and negative ranges respectively. Moreover, as shown in
\cite{MR588409}, this discussion can be passed to a continuous time limit to derive the companion circle
of fluctuation identities for maxima, minima  and ladder processes associated with L\'evy processes. 
In section \ref{sec:max} we give new explanations for the compound Poisson distributions mentioned above.

The rest of this article is structured as follows. In Section \ref{sec:proof} we will prove Theorem \ref{theorem:1} under assumption $\mathbf{A}$ and give corollaries relating to the partition and composition induced by the concave majorant. In Section \ref{sec:examples} we will analyze some specific examples of composition probabilities, including the Cauchy increment case, which turns out to be particularly simple. In Section \ref{sec:poisson} we extend the description to the case where $n$ is replaced by $n(q)$, a geometric random variable with parameter $1-q$, which results in a description of the concave majorant and the excursions under each face as a Poisson point process. In Section \ref{sec:applications} we apply the Poissonian theory. First, by letting $q \rightarrow 1$ we find a description of the concave majorant for the random walk on $[0, \infty)$, and the associated excursions under each face. Then we analyze the behaviour of the concave majorant as $n$ grows. As a final application we investigate the pre and post maximum parts of the walk. In Section \ref{sec:max} we investigate the two concave majorants that result from decomposing the random walk at its maximum, and their associated partitions. In Section \ref{sec:noncontinuous} we extend the theory to $X_1,\ldots,X_n$ not satisfying assumption $\mathbf{A}$. Also in Section \ref{sec:noncontinuous} we investigate three nested compositions of integers that arise naturally. At the end of this  Section \ref{sec:noncontinuous}  some examples of how the general theory can be applied are given. In Section \ref{sec:conditionedwalk} we finish answering the rearrangement problem mentioned above by describing the law of a random walk conditional on the value of its concave majorant. Finally, in Section \ref{sec:pathtransform}, we describe an important path transformation that provides Pitman and Uribe Bravo with the basis for a full investigation into the concave majorant of a L\'evy process \cite{pitmanbravo}.

\section{Proof of Theorem \ref{theorem:1} under assumption $\mathbf{A}$ and the partition and composition laws}
\label{sec:proof}

We begin with a simple Lemma due to Spitzer relating to cyclic permutations of increments of walks that shows that under assumption $\mathbf{A}$ the appropriate cyclic permutations discussed in the introduction are almost surely unique.

\begin{lemma}\cite[Theorem $2.1$]{MR0079851}
\label{lem:cyclics}
Let $x = (x_1, \ldots , x_n)$ be a vector such that no two subsets of the coordinates have the same arithmetic mean. For  $1 \le k \le n$ let $x_{k+n} = x_k$, and let $x(k) =(x_k, x_{k+1}, \ldots , x_{k+n})$. Then there is a unique $1 \le k^* \le n$ such that the walk with increments $x(k^*) =(x_{k^*}, x_{k^*+1}, \ldots , x_{k^*+n})$ lies below the chord  joining its start and end points.
\end{lemma}

\noindent{\textbf{Proof. (Theorem \ref{theorem:1} under assumption $\mathbf{A}$)} By conditioning on the set of values that $X_1,\ldots,X_n$ take it is enough to show that $S_{\rho}^{[0,n]} \ed S^{[0,n]}$ in the case where $X_1,\ldots,X_n$ are samples without replacement from $n$ real numbers $x_1,\ldots,x_n$ such that no two subsets of $x_1,\ldots,x_n$ have the same arithmetic mean. Thus it is enough to show that for every permutation $\sigma \in \Sigma_n$ we have
$$
\P(X_{\rho(1)} = x_{\sigma(1)}, \ldots , X_{\rho(n)} = x_{\sigma(n)} ) = \frac{1}{n!}
$$
and without loss of generality it is enough to show this for $\sigma$ the identity permutation.
Suppose the concave majorant  of the deterministic walk with increments $(x_1, \ldots , x_n ) $ has $k$ faces whose lengths \emph{in order of appearance} are $(m_1,\ldots,m_k)$, so that the composition induced by the lengths of the faces of the concave majorant is $(m_1, \ldots,m_k)$. Let $\tau \in \Sigma_k$ be such that
$$
(n_1,\ldots,n_k) := (m_{\tau(1)},\ldots , m_{\tau(k)}) 
$$
are the lengths of the $k$ faces in \emph{non-increasing order}, so that the partition induced by the lengths of the faces of the concave majorant is $(n_1,\ldots,n_k)$.

First suppose that each element of $(n_1,\ldots,n_k)$ is distinct. Then the event $\{ X_{\rho( \ell )} = x_{\ell} \: : \: 1 \le \ell \le n \}$ occurs if and only if
\begin{enumerate}[(i)]
\item the partition chosen according to the lengths of the cycles of a random permutation is $(n_1,\ldots,n_k)$;
\item for each $1 \le i \le k$, the ordered list $ (X_{n_1+\cdots+n_{i-1}+1}, \ldots, X_{n_1+\cdots+n_i})$ is one of the $n_i$ cyclic permutations of the ordered list 
\newline $( x_{m_1+m_2+\cdots+m_{\tau(i)-1}+1}, \ldots , x_{m_1+m_2+\cdots+m_{\tau(i)}} ) $.
\end{enumerate}
According to the Ewens Sampling Formula, the event in (i) has probability $\prod_{i=1}^k \frac{1}{n_i}$. The event in (ii) is independent 
of the event in (i), and has probability $\frac{1}{n!} \prod_{i=1}^k n_i$.

Now suppose that the elements of $(n_1,\ldots,n_k)$ are not distinct. For $1 \le j \le n$ let $I_j = \{ i \: : \: n_i = j \}$ and let $a_j = |I_j|$.
The event $\{ X_{\rho( \ell )} = x_{\ell} \: : \: 1 \le \ell \le n \}$ 
occurs if and only if
\begin{enumerate}[(i)]
\item the partition chosen according to the lengths of the cycles of a random permutation is $(n_1,\ldots,n_k)$;
\item for each $1 \le j \le n$, for each $i \in I_j$ the ordered list $(X_{n_1+\cdots+n_{i-1}+1}, \ldots, X_{n_1+\cdots+n_i})$ is one of the $n_i=j$ cyclic permutations of the ordered list 
\newline $( x_{m_1+m_2+\cdots+m_{\tau(i')-1}+1}, \ldots , x_{m_1+m_2+\cdots+m_{\tau(i')}} ) $ for some $i' \in I_j$.
\end{enumerate} 
By the Ewens Sampling Formula, the event in (i) has probability
$
\left( \prod_{i=1}^k \frac{1}{n_i} \right) \left( \prod_{j=1}^n \frac{1}{a_j!} \right).
$
The event in (ii) is independent of the event in (i), and has probability 
$
\frac{1}{n!} \left( \prod_{i=1}^k n_i \right) \left( \prod_{j=1}^n a_j ! \right).
$
Hence $\P ( X_{\rho( \ell )} = x_{\ell} \: : \: 1 \le \ell \le n ) = \frac{1}{n!}$.
\endpf
\vspace{10pt}

As a direct consequence of Theorem \ref{theorem:1} we have the result of Goldie \cite{MR994088} mentioned in the introduction. 

\begin{corollary}\label{cor:goldie}
Let $M_{n,1}, \ldots, M_{n,F_n}$ be the lengths of the faces of the concave majorant
of $S^{[0,n]}$ arranged in non-increasing order.
Then under assumption $\mathbf{A}$ the joint distribution of $M_{n,1}, \ldots, M_{n,F_n}$ is given
by the formula 
$$ \P( F_n = k, M_{n,i} = n_i, 1 \le i \le k) = \prod_{j=1}^n \frac{1}{j^{a_j} a_j !} $$
for all $(n_1,\ldots,n_k) \in \mathcal{P}_n$, where $a_j = \# \{ i \: : \: 1 \le i \le k , n_i = j \}$ for $1 \le j \le n$. I.e.\ The partition of $n$ induced by the lengths of the faces of the concave majorant of $S^{[0,n]}$ has the law of a partition of $n$ induced by the cycle lengths of a random permutation.
\end{corollary}
\proof
Following the construction in Theorem \ref{theorem:1}, the lengths $L_{n,1},\ldots,L_{n,K_n}$ are exactly the lengths of the faces of the concave majorant of $S_{\rho}^{[0,n]}$, and the conclusion follows since $S^{[0,n]} \ed S_{\rho}^{[0,n]}$.
\endpf
\vspace{10pt}

Further, Theorem \ref{theorem:1} allows us to describe the law of the composition induced by the lengths of the faces of the concave majorant.

\begin{corollary}\label{cor:composition}
Let $(N_{n,1}, \ldots, N_{n,F_n})$ be the composition of $n$ induced by the lengths  of the faces of the concave majorant of $S^{[0,n]}$.
Then under assumption $\mathbf{A}$ 
the joint distribution of $N_{n,1}, \ldots, N_{n,F_n}$ is given
by the formula
$$
\P( F_n = k, N_{n,i} = n_i, 1 \le i \le k) = 
\P \left( \frac { S_{n_1} ^{(1)} } { n_1 }  > \frac { S_{n_2} ^{(2)} } { n_2 }  >  \cdots > \frac { S_{n_k} ^{(k)}} { n_k } \right) \prod_{i = 1}^k \frac{1}{n_i}
$$
for all compositions $(n_1, \ldots, n_k)$ of $[n]$ into $k$ parts,
where for $1 \le i \le k$
$$
S_{n_i} ^{(i)}:=  S_{n_1 + \cdots + n_{i} } - S_{n_1 + \cdots + n_{i-1} } \ed S_{n_i}
$$
In particular, if the $X_i$ are independent, then so are the $S_{n_i}^{(i)}$ for $1 \le i \le k$.
\end{corollary}
\proof
Fix a composition $(n_1,\ldots,n_k)$ and let  
$(\overrightarrow{n}_{\tau(1)} ,\ldots,\overrightarrow{n}_{\tau(k)} )$ be $(n_1,\ldots,n_k)$ in non-increasing order.
Let $T$ be the set of $\tau \in \Sigma_k$ such that $(\overrightarrow{n}_{\tau(1)},\ldots,\overrightarrow{n}_{\tau(k)}) = (n_1,\ldots,n_k)$. 
Then $|T| = \prod_{j=1}^n a_j$, where $a_j = \# \{ i \: : \: 1 \le i \le k , n_i = j \}$ for $1 \le j \le n$.
We are interested in comparing the slopes of the faces of the concave majorant that result from the construction in Theorem \ref{theorem:1}. In this direction, for $1 \le i \le k$ let
$$
S_{\overrightarrow{n}_{\tau(i)}}^{(\tau(i))} =
S_{\overrightarrow{n}_1 + \cdots + \overrightarrow{n}_{\tau(i)}} - S_{\overrightarrow{n}_1 + \cdots + \overrightarrow{n}_{\tau(i)-1}} \ed S_{\overrightarrow{n}_{\tau(i)}} = S_{n_i}
$$
Under the construction in Theorem \ref{theorem:1}, the events $\{ F_n = k \}$ and $ \{ N_{n,i} = n_i \: : \: 1 \le i \le k \} $ occur if and only if
\begin{enumerate}[(i)]
\item $(L_{n,1} , \ldots, L_{n,K_n}) = (\overrightarrow{n}_1,\ldots,\overrightarrow{n}_k)$;
\item  $\frac { S_{ \overrightarrow{n}_{\tau(1)} }^{(\tau(1))} } { n_1 }  > \frac { S_{ \overrightarrow{n}_{\tau(2)} }^{(\tau(2))} } { n_2 } >  \cdots > \frac { S_{ \overrightarrow{n}_{\tau(k)} }^{(\tau(k))} } { n_k }$ for some $\tau \in T$.
\end{enumerate}
As before, the event in (i) has probability
$
\left( \prod_{i=1}^k \frac{1}{n_i} \right) \left( \prod_{j=1}^n \frac{1}{a_j!} \right) .
$
The event in (ii) is independent of the event in (i), and by exchangeability the probability that it occurs for one particular element of $T$ is 
$$
\P \left( \frac { S_{n_1} ^{(1)} } { n_1 }  > \frac { S_{n_2} ^{(2)} } { n_2 }  >  \cdots > \frac { S_{n_k} ^{(k)}} { n_k } \right)  
$$
Recalling that $|T| = \prod_{j=1}^n a_j$ completes the proof.
\endpf

\section{Examples of composition probabilities}
\label{sec:examples}

The special case of Cauchy increments gives rise to the following appealing version of Corollary \ref{cor:composition}.

\begin{corollary}
\label{cor:cauchy}
Suppose that the $X_i$ are independent and such that
$S_k/k$ has the same distribution for every $k$, as when
the $X_i$ have a Cauchy distribution. Then
$$
\P( F_n = k; N_{n,i} = n_i, 1 \le i \le k) = 
\frac{1 }{k!} \prod_{i = 1}^k \frac{1}{n_i}
$$
and hence $\{ N_{n,i} \: : \: 1 \le i \le F_n \}$ has the same distribution as the composition of $n$ created by first choosing a random permutation of $n$ and then putting the cycle lengths in uniform random order.
\end{corollary}
\proof
Since $\frac { S_{n_1} ^{(1)} } { n_1 }  ,  \ldots , \frac { S_{n_k} ^{(k)}} { n_k } $ is an i.i.d. sequence each of the $k!$ orderings is equally likely, and hence $\P(\frac { S_{n_1} ^{(1)} } { n_1 }  > \cdots > \frac { S_{n_k} ^{(k)}} { n_k } ) = \frac{1}{k!}$.
\endpf

\vspace{10pt}

Note that the continuum limit of this result can be read from Bertoin's work
\cite{MR1747095}. 
The above result shows that the
Cauchy discrete model is the same as that derived by random sampling from the continuum Cauchy model, as per Gnedin's
theory of sampling consistent compositions of positive integers \cite{gp03}. 
That is, let $U_1,\ldots,U_n$ be independent identically distributed uniform random variables on $[0,1]$ and let $X$ be a Cauchy process on $[0,1]$. Generate a composition of $n$ by putting $i$ in the same block as $j$ if and only if $U_i$ and $U_j$ fall in the same segment of the composition of $[0,1]$ induced by the lengths of the faces of the concave majorant of $X$, and then ordering blocks according to the ordering of the faces of the concave majorant of $X$. Then the composition of $n$ that is generated will have the same distribution as $(N_{n,1},\ldots,N_{n,F_n})$ in Corollary \ref{cor:cauchy}.
This does not seem at all obvious a priori,
and according to simulation is not true in the Brownian case, suggesting that it is not true in general.

Now let $X_1,\ldots,X_n$ be any exchangeable sequence of random variables satisfying assumption $\mathbf{A}$, as in Corollary \ref{cor:composition}. We now give some numerical examples of composition probabilities when $n$ is small. Let
$$
p(n_1, \ldots, n_k):= \P( F_n = k, N_{n,i} = n_i, 1 \le i \le k) 
$$
Using symmetry and the partition probabilities given in Corollary \ref{cor:goldie}, universal values are
$$
p(1,1) = 1/2, \: \: \: p(2)  = 1/2
$$
$$
p(3) = 1/3, \: \: \: p(2,1) = p(1,2) = 1/4, \: \: \: p(1,1,1) = 1/6
$$
$$
p(4) = 1/4,  \: \: \: p(1,3) = p(3,1) = 1/6, \: \: \: p(2,2) = 1/8, \: \: \: p(1,1,1,1) = 1/24
$$
As $n$ increases, the first values that depend on the particular choice of increment distributions are
$$
p(1,1,2) =  p(2,1,1) = \frac{1}{2} \P( X_1 > X_2 > \mbox{$\frac{1}{2}$}(X_3 + X_4) ) 
$$
$$
p(1,2,1) = \frac{1}{2} \P( X_1 > \mbox{$\frac{1}{2}$}(X_2 + X_3 ) > X_4 )
$$
where according to the partition probabilities we must have
$$
p(1,1,2) +  p(2,1,1) + p(1,2,1) = 1/4
$$
We consider two special cases - independent Cauchy increments and independent Gaussian increments. When the increments are independent and Cauchy, the 3 probabilities above are equal, with
$$
2p(1,2,1) = \P( X_1 > \mbox{$\frac{1}{2}$}(X_2 + X_3 ) > X_4 ) = 1/6 = 0.1666666...
$$
Note that
$$
\P( X_1 > \mbox{$\frac{1}{2}$}(X_2 + X_3 ) > X_4 ) =  \P( \mbox{$\frac{1}{2}$}(X_2 + X_3 ) - X_1 < 0 \mbox { and } X_4 - \mbox{$\frac{1}{2}$}(X_2 +X_3) < 0  ).
$$
In the centered Gaussian case with $Var(X_1) = 1$ this is the probability of the negative quadrant for a centered bivariate normal with equal variances $3/2$
and covariance $-1/2$ and thus correlation $\rho = -1/3$. That probability is given by
$$
\frac{1 }{4 } + \frac{ \arcsin(-1/3) }{ 2 \pi } = 0.195913276
$$
The difference with the Cauchy case is quite small. The fact that it is larger
is consistent with the known differences in behaviour of the limit partitions for large $n$ after scaling; it is known that the concave majorant of Brownian motion is more likely to have longer faces in its central region than the concave majorant of a Cauchy process. We conclude this section by conjecturing that $p(1,2,1)$ is a monotonic function of the stability index $\alpha$ for symmetric stable laws.

\section{A Poisson point process description}
\label{sec:poisson}

The concave majorant of $S^{[0,n]}$ can be viewed as a random point process on $\{ 1, \ldots , n \} \times \mathbb{R}$, where a point at $(j,s)$ means that one of the faces of the concave majorant has length $j$ and increment $s$.
Let $A_n(j)$ be the number of faces of the concave majorant of $S^{[0,n]}$ that have length $j$ for $1 \le j \le n$, and let $\Sigma_j^{(1)}, \ldots , \Sigma_j^{(A_n(j))}$ be the increments of the faces with length $j$ in uniform random order. Thus if $X_1, \ldots , X_n$ are independent then for each $1 \le j \le n$, conditionally given $A_n(j) = a_j$, $\Sigma_j^{(\ell)}$ is an independent copy of $S_j$ for each $1 \le \ell \le a_j$. Figure \ref{fig:pointprocess} shows an example of such a point process.
To construct the concave majorant from this point process the faces with lengths and increments indicated by the points are arranged in decreasing order of slope.

\begin{figure}
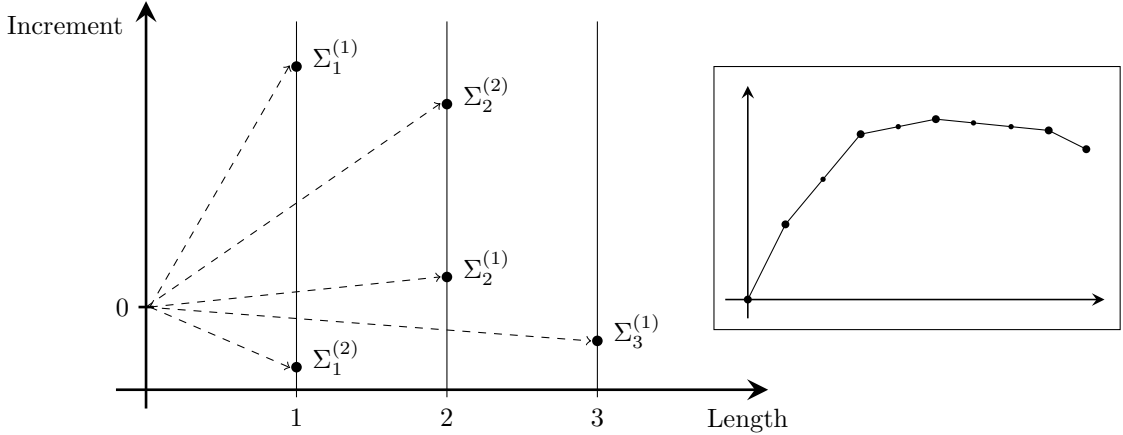


\begin{pgfpicture}{0cm}{0cm}{20.5cm}{6.5cm}
\begin{pgfscope} 
\pgfsetendarrow{\pgfarrowsingle}
\pgfsetlinewidth{1pt} 
\pgfxyline(-0.4,0.5)(8.2,0.5) 
\pgfxyline(-0,0.25)(-0,5.6) 
\end{pgfscope} 
\begin{pgfscope} 
\pgfsetlinewidth{1pt} 
\pgfxyline(-0.1,1.6)(0.1,1.6) 
\end{pgfscope} 
\pgfline{\pgfxy(2,0.4)}{\pgfxy(2,5.4)} 
\pgfline{\pgfxy(4,0.4)}{\pgfxy(4,5.4)} 
\pgfline{\pgfxy(6,0.4)}{\pgfxy(6,5.4)} 
\pgfcircle[fill]{\pgfpoint{2cm}{4.8cm}}{2pt}
\pgfcircle[fill]{\pgfpoint{2cm}{0.8cm}}{2pt}
\pgfcircle[fill]{\pgfpoint{4cm}{2cm}}{2pt}
\pgfcircle[fill]{\pgfpoint{4cm}{4.3cm}}{2pt}
\pgfcircle[fill]{\pgfpoint{6cm}{1.15cm}}{2pt}
\pgfputat{\pgfxy(2,0.25)}{\pgfbox[center,top]{1}}
\pgfputat{\pgfxy(4,0.25)}{\pgfbox[center,top]{2}}
\pgfputat{\pgfxy(6,0.25)}{\pgfbox[center,top]{3}}
\pgfputat{\pgfxy(8,0)}{\pgfbox[center,base]{Length}}
\pgfputat{\pgfxy(-0.3,5.35)}{\pgfbox[right,center]{Increment}}
\pgfputat{\pgfxy(-0.4,1.6)}{\pgfbox[fight,center]{0}}
\pgfputat{\pgfxy(2.1,4.7)}{\pgfbox[left,bottom]{  $\Sigma_1^{(1)}$  }}
\pgfputat{\pgfxy(2.1,0.7)}{\pgfbox[left,bottom]{  $\Sigma_1^{(2)}$  }}
\pgfputat{\pgfxy(4.1,4.2)}{\pgfbox[left,bottom]{  $\Sigma_2^{(2)}$  }}
\pgfputat{\pgfxy(4.1,1.9)}{\pgfbox[left,bottom]{  $\Sigma_2^{(1)}$  }}
\pgfputat{\pgfxy(6.1,1.05)}{\pgfbox[left,bottom]{  $\Sigma_3^{(1)}$  }}
\begin{pgfscope}
\pgfsetendarrow{\pgfarrowto}
\pgfsetdash{{3pt}{3pt}}{0pt}
\pgfxyline(0.05,1.6)(1.9,4.8)
\pgfxyline(0.05,1.6)(1.9,0.8) 
\pgfxyline(0.05,1.6)(3.9,2) 
\pgfxyline(0.0,1.6)(3.9,4.3)
\pgfxyline(0.0,1.6)(5.9,1.15)
\end{pgfscope}
\begin{pgfscope} 
\pgfsetendarrow{\pgfarrowsingle}
\pgfsetlinewidth{0.6pt} 
\pgfxyline(7.7,1.7)(12.7,1.7) 
\pgfxyline(8,1.45)(8,4.5) 
\end{pgfscope} 
\pgfrect[stroke]{\pgfpoint{7.55cm}{1.3cm}}{\pgfpoint{5.40cm}{3.5cm}}
\pgfxyline( 8 , 1.7 )( 8.5 , 2.7 ) 
\pgfxyline( 8.5 , 2.7 )( 9 , 3.3 ) 
\pgfxyline( 9 , 3.3 )( 9.5 , 3.9 ) 
\pgfxyline( 9.5 , 3.9 )( 10 , 4 ) 
\pgfxyline( 10 , 4 )( 10.5 , 4.1 ) 
\pgfxyline( 10.5 , 4.1 )( 11 , 4.05 ) 
\pgfxyline( 11 , 4.05 )( 11.5 , 4 ) 
\pgfxyline( 11.5 , 4 )( 12 , 3.95 ) 
\pgfxyline( 12 , 3.95 )( 12.5 , 3.7 )
\pgfcircle[fill]{\pgfpoint{8cm}{1.7cm}}{1.5pt}
\pgfcircle[fill]{\pgfpoint{8.5cm}{2.7cm}}{1.5pt}
\pgfcircle[fill]{\pgfpoint{9.5cm}{3.9cm}}{1.5pt}
\pgfcircle[fill]{\pgfpoint{10.5cm}{4.1cm}}{1.5pt}
\pgfcircle[fill]{\pgfpoint{12cm}{3.95cm}}{1.5pt}
\pgfcircle[fill]{\pgfpoint{12.5cm}{3.7cm}}{1.5pt}
\pgfcircle[fill]{\pgfpoint{9cm}{3.3cm}}{1pt}
\pgfcircle[fill]{\pgfpoint{10cm}{4cm}}{1pt}
\pgfcircle[fill]{\pgfpoint{11cm}{4.05cm}}{1pt}
\pgfcircle[fill]{\pgfpoint{11.5cm}{4cm}}{1pt}
\end{pgfpicture}

\caption{An example point process and the resulting concave majorant. The dashed lines show the slope of each face, and these faces are arranged in decreasing order of slope.}
\label{fig:pointprocess}
\end{figure}

Now suppose we have an infinite sequence of exchangeable random variables $X_1, X_2, \ldots $, such that almost surely no two subsets have the same arithmetic mean. As before let $S_0=0$ and $S_j = \sum_{i=1}^j X_i$ for $j \ge 1 $. Following ideas from the fluctuation theory of Greenwood and Pitman \cite{greenpit80} we now randomise the length of the walk by setting the number of steps of the random walk equal to $n(q)$, where $n(q)$ is a geometric random variable with parameter $1-q$, so that
$$\P(n(q) \geq n) = q^n \quad \text{ for }n=0,1,2,\ldots$$
Let $S^{[0,n(q)]} = \{ (j,S_j) \: : \: 0 \le j \le n(q) \}$, and let
$$
0 < N_{n(q),1} < N_{n(q),1} + N_{n(q),2} < \cdots < N_{n(q),1} + \cdots + N_{n(q),F_{n(q)}} = n(q)
$$
be the successive times that $S^{[0,n(q)]}$ meets its concave majorant, where $F_{n(q)}$ is the number of faces of the concave majorant of $S^{[0,n(q)]}$. The following Lemma, which involves a fundamental Poisson representation of the geometric distribution, is due to Shepp and Lloyd \cite{MR0195117}, who were just working with partitions generated by random permutations, not concave majorants.

\begin{lemma}\label{lem:poisson1}
Let $A_j = \# \{ i \: : \: 1 \le i \le F_{n(q)} , N_{n(q),i} = j \}$ for $j \ge 1$. Then $A_j$ has the Poisson distribution with mean $q^j/j$, independently for each $j \ge 1$.
\end{lemma}
\proof
Noting that $\log(1-q) = - \sum_j q^j/j$, we have that
\begin{eqnarray*}
\P(A_j = a_j , j \ge 1) & = & \P( n(q) = \mbox{$\sum_{j \ge 1}$} ja_j ) \P ( A_j = a_j , j \ge 1 | n(q) = \mbox{$\sum_{j \ge 1}$} j a_j ) \\
& = & (1-q) q^{\sum_j j a_j} \frac{1}{\prod_j j^{a_j} a_j ! } \\
& = & \prod_j \frac{ \left( \frac{q^j }{j} \right)^{a_j} e^{-\frac{q^j}{j}} }{ a_j!} 
\end{eqnarray*}
where the second equality comes from Corollary \ref{cor:goldie}.
\endpf
\vspace{10pt}

For the next theorem, and in fact the rest of this section, it is important that we assume $X_1, X_2, \ldots $ are independent with common continuous distribution. The theorem asserts that the point process discussed above is a Poisson point process under this assumption.

\begin{theorem}\label{theorem:poisson1}
If $X_1, X_2, \ldots $ are independent with common continuous
distribution, then the point process of lengths and increments of
faces of the concave majorant of $S^{[0,n(q)]}$ is a Poisson point process on $\{ 1,2,\ldots \} \times \mathbb{R}$ with intensity $j^{-1} q^j \P(S_j \in dx)$ for $j = 1, 2, \ldots$, $x \in \mathbb{R}$. Moreover, let 
$T_i =  \sum_{l = 1}^i N_{n(q),l}$, $0 \le i \le F_{n(q)}$, be the consecutive times at which $S^{[0,n(q)]}$ meets its concave majorant, so that $T_0=0$ and $T_{F_{n(q)}} = n(q)$. Then 
the sequence of path segments
$$
\{ (S_{T_i + k} - S_{T_i},  0 \le k \le N_{n(q),i} ) , i = 0, \ldots , F_{n(q)} -1 \},
$$
is a list of the points of a Poisson point process in the space of finite random walk segments
$$
\{ (s_1, \ldots , s_j) \mbox{ for some } j = 1,2, \ldots \}
$$
whose intensity measure on paths of length $j$ is $q^j j^{-1}$ times the conditional distribution of $(S_1, \ldots , S_j)$ given that $S_k < (k/j) S_j $ for all $1 \le k \le j-1 $.
\end{theorem}
\proof
Conditionally given $A_j = a_j$ the increment for each face of length $j$ is an independent copy of $S_j$ by Theorem \ref{theorem:1}. Combined with Lemma \ref{lem:poisson1} this proves the first statement.

Conditional on the concave majorant of $S^{[0,n(q)]}$ having a face of length $j$ and increment $s$, the increments of $S^{[0,n(q)]}$ over that face of the concave majorant have the distribution of $(X_1, \ldots , X_j)$ given that $\sum_{\ell=1}^k X_{\ell} < (k/j) s$ for all $1 \le k \le j-1 $ and $\sum_{\ell=1}^j X_{\ell} = s$, and this law is independent for each face of $S^{[0,n(q)]}$. This implies the second statement.
\endpf
\vspace{10pt}

A simple but important corollary of Theorem \ref{theorem:poisson1} is the following.

\begin{corollary}
\label{cor:compoundpois}
$(n(q),S_{n(q)})$ has a compound Poisson distribution, and  the total number of faces 
$F_{n(q)}$ of the concave majorant of $S^{[0,n(q)]}$ has Poisson distribution with mean
$$
\sum_{j = 1}^\infty j^{-1} q^j= -\log(1-q).
$$
\end{corollary}

The first assertion of Corollary \ref{cor:compoundpois} can in fact be
seen directly since $(n(q),S_{n(q)}) = \sum_{i=1}^{n(q)} (1,X_i)$ and
$n(q)$ is itself compound Poisson. Explicitly, $n(q)$ is a Poisson
compound of a log-series law: $n(q)$ has probability generating
function $\E n(q) = (1 - q)/(1 - qz)$ which can be expressed as 
$e^{- \lambda(1-\phi(z))}$ where $\lambda = -\ln (1 - q)$ and $\phi$ is
  the probability generating function of the log-series law with
  parameter $q$. This well known decomposition of a geometric random
  variable reappears later in Lemma \ref{lem:E_qj}.

\section{Applications of the Poissonian description}
\label{sec:applications}

\subsection{The random walk on $[0,\infty)$}
\label{sec:infwalk}

By letting $q \rightarrow 1$ it is possible to deduce the structure of the concave majorant of the random walk on $[0, \infty)$ using Theorem \ref{theorem:poisson1}. Groeneboom \cite{MR714964} gave a Poissonian description of the concave majorant of BM on $[0,\infty)$; that there is a closely parallel description for random walks does not seem to have been pointed out before. The case of L\'evy processes will be covered in the forthcoming paper by Pitman and Uribe Bravo \cite{pitmanbravo}.

Suppose $\mathbb{E}(X_1) = \mu \in [-\infty,\infty)$. Informally, as $q \rightarrow 1$ the intensity measure of the Poisson point process of face lengths and increments approaches $j^{-1} \P(S_j \in dx)$, but since the slope of the concave majorant converges downwards to $\mu$ but does not reach it, only the faces with slope greater than $\mu$ will contribute to the concave majorant in the limit. Therefore by Poisson thinning we get a new intensity measure $j^{-1} \P(S_j \in dx) 1(x > j \mu)$. Moreover, we can also describe path segments of the walk below each face of the concave majorant as a Poisson point process.

\begin{theorem}
\label{theorem:poisson2}
Let $S_0=0$ and $S_j = \sum_{i=1}^j X_i$ for $j \ge 1 $, where $X_1, X_2, \ldots $ are independent random variables with common continuous distribution that has a well defined mean $\mu:= \E(X_1) \in [-\infty, \infty)$. Let $S^{[0,\infty)} = \{ (j,S_j) \: : \: j \ge 0 \} $. Let $0  = T_0 < T_1 < T_2 < \cdots$ be the successive times that $S^{[0,\infty)}$ meets its concave majorant, and let $N_i = T_i-T_{i-1}$ for $i \ge 1$. 
Then the sequence of path segments
$$
\{ (S_{T _i + k} - S_{T_i},  0 \le k \le N_i ) , i = 0,2, \ldots \}
$$
is a list of the points of a Poisson point process in the space of finite random walk segments
$$
\{ (s_1, \ldots , s_j) \mbox{ for some } j = 1,2, \ldots \}
$$
whose intensity measure on paths of length $j$ is $j^{-1}$ times the restriction to $S_j \in (j \mu,\infty)$
of the conditional distribution of $(S_1, \ldots , S_j)$ given that $S_k < (k/j) S_j $ for all $1 \le k < j $.
\end{theorem}
\proof
The combination of the following four facts is enough to prove the theorem:
\begin{enumerate}[(i)] 
\item the number of faces of length $j$ has a Poisson distribution with mean $j^{-1} \P ( S_j > j \mu )$;
\item these numbers are independent as $j$ varies;
\item given all of these numbers, and with $n$ faces of length $j$, the $n$ walks on the associated faces, when listed in a uniform random order independently of the walks on the faces, are $n$ independent processes each distributed according to $(S_1, \ldots , S_j)$ given that $S_k < (k/j) S_j $ for all $1 \le k < j $ and $S_j > j \mu$.
\item given $n$ faces of length $j$, the increments of these faces, when listed in uniform random order, are distributed like $n$ independent copies of $S_j$ given $S_j > j \mu$.
\end{enumerate}
The main thing to check is that (i) and (ii) are true, i.e.\ that the counts
$$
A_\infty  (j) := \# \{ j : N_i = j \} 
$$
are independent Poisson variables with mean $j^{-1} \P ( S_j \ge j \mu )$.
Once we have shown this, (iii) and (iv) follow from Poisson thinning and previous discussions relating to the independence of the walks below each segment.

Let $n(q)$ be a geometric random variable with parameter $1-q$. Let $S^{[0,n(q)]} = \{ (j,S_j) \: : \: 0 \le j \le n(q) \}$, so that the concave majorant of $S^{[0,n(q)]}$ and $S^{[0,\infty)}$ agree up until some random time $T_{n(q)}^*$.
\begin{lemma}
\label{lem:lastvertex}
$T_{n(q)}^*$ is the maximal $T_i$ with $T_i \le n(q)$.
\end{lemma}
\proof
To see this, let $i$ be such that $T_i \le n(q)$. Since the concave
majorant of $S^{[0,n(q)]}$ is everywhere less than or equal to the concave majorant of $S^{[0,\infty)}$, if they did not agree at time $T_i$ then the concave majorant of $S^{[0,n(q))}$ would go beneath the point $(T_i,S_{T_i})$, but this is a contradiction since $(T_i,S_{T_i})$ is in $S^{[0,n(q)]}$.
\endpf
\vspace{10pt}

Let
$$
A_{n(q)}(j) := \# \{ i : N_{n(q),i} = j \} 
$$
where  $N_{n(q),1}, \ldots , N_{n(q),F_{n(q)}}$ are the lengths of faces of the concave majorant of $S^{[0,n(q)]}$. There are the obvious decompositions
\begin{eqnarray}
A_{\infty}(j) & = & A_{\infty}(j) (0, T_{n(q)}^*] + A_{\infty}(j) (T_{n(q)}^*,\infty] \label{eq:decomp1} \\
A_{n(q)}(j) & = & A_{n(q)}(j) (0, T_{n(q)}^*] + A_{n(q)}(j) (T_{n(q)}^*,\infty] \label{eq:decomp2}
\end{eqnarray}
where  e.g. $A_{\infty}(j) (0, T_{n(q)}^*]$ is the number of faces of the concave majorant of $S^{[0,\infty)}$  of length $j$ up to  and including the face ending at time $T_{n(q)}^*$, and the other terms are defined similarly. Moreover, 
since $T_{n(q)}^*$ is by definition the maximal common vertex of the concave majorants of $S^{[0,n(q)]}$ and $S^{[0,\infty)}$, it is clear that 
\begin{eqnarray}
A_{\infty}(j) (0, T_{n(q)}^*]  \! \! & \! = \! & \! A_{n(q)}(j) (0, T_{n(q)}^*] \nonumber \\
& \! = \! & \! \# \{ i : N_{n(q),i} = j, \: S_{T_i}- S_{T_{i-1}} > j \alpha_{n(q)} \} \label{eq:poisson} 
\end{eqnarray}
where $\alpha_{n(q)}$ is the right derivative of the concave majorant of $S^{[0,\infty)}$ at time $T_{n(q)}^*$. Conditionally given $\alpha_{n(q)}$, by Poisson thinning and Theorem \ref{theorem:poisson1} the distribution of the right hand side of \re{eq:poisson} is Poisson with mean $q^jj^{-1} \P(S_j > j \alpha_{n(q)})$, independently for each $j$. The strategy at this point is to let $q \te 1$, so that $T_{n(q)} \te \infty$ and $\alpha_{n(q)} \te \mu$, resulting in $A_{\infty}(j) $ having Poisson distribution with mean $j^{-1} \P( S_j > j \mu)$, independently for each $j$, i.e.\ resulting in (i) and (ii).

Let $\{q_m \}_{m \geq 1}$ be any sequence such that if $\{ n(q_m) \}_{m \ge 1}$ is a sequence of independent geometric random variables with parameters $1-q_m$ then $n(q_m) \te \infty$ almost surely as  $m \te \infty$ (so that necessarily $q_m \te 1$). Suppose that  $T_{(n(q_m))} \te \infty$ and $\alpha_{n(q_m)} \te \mu$ almost surely, so that
\begin{eqnarray}
A_{\infty}(j)   &  =  & \lim_{m \te \infty} A_{\infty}(j) (0, T_{(n(q_m))}] \nonumber \\
&  =  & \lim_{m \te \infty} \# \{ i : N_{n(q_m),i} = j, \: S_{T_i}- S_{T_{i-1}} > j \alpha_{n(q_m)} \} \label{eq:endpois} \\
\end{eqnarray}
where the first equality is from \re{eq:decomp1} and the second is from \re{eq:poisson}. Since $\alpha_{n(q_m)} \te \mu$ almost surely, by continuity of the function $x \mapsto \P(S_j > jx)$ the distribution of the right hand side of \re{eq:endpois} is Poisson with parameter $j^{-1} \P(S_j > j \mu)$, independently for each $j$. This proves (i) and (ii).

It remains to prove that  $T_{(n(q_m))} \te \infty$ and $\alpha_{n(q_m)} \te \mu$ almost surely as $m \te \infty$. For every $i \geq 1$, since $T_i < \infty$ we will have $n(q_m) > T_i$ eventually, and hence by Lemma \ref{lem:lastvertex}  for every $i \geq 1$ we will have $T_{(n(q_m))} \geq T_i$ eventually. Since $T_i \te \infty$ this implies that $T_{(n(q_m))} \te \infty$ almost surely.

\begin{lemma}
\label{lem:lessthanmu}
Almost surely no face of the concave majorant of $S^{[0,\infty)}$ can have slope less than $\mu$.
\end{lemma}
\proof
If $\mu = - \infty$ then the conclusion is clear. Suppose $\mu \in (-\infty,\infty)$, then since $S_n - n \mu$ is a mean zero random walk and hence recurrent, for every $i \geq 1$ there will almost surely be some $n_i > T_i$ such that $S_{n_i} > S_{T_i} + (n_i-T_i) \mu$, and hence for any vertex of the concave majorant the slope of the face to the right must be greater than $\mu$.
\endpf
\vspace{10pt}

\begin{lemma}
\label{lem:morethanmu}
For every $\epsilon > 0$ there will almost surely be a face of the concave majorant with slope $x$ such that $\mu < x < \mu + \epsilon$.
\end{lemma}
\proof
For any $\mu \in [-\infty,\infty)$ by the strong law of large numbers $S_n/n \te \mu$ almost surely as $n \te \infty$. But if there was no slope of the concave majorant on $[0, \infty)$ with slope $x < \mu + \epsilon$ then we would have $\limsup_n S_n/n > \mu$. Combined with Lemma \ref{lem:lessthanmu} this gives the conclusion.
\endpf
\vspace{10pt}

We already have that $T_{(n(q_m))} \te \infty$ almost surely. Since
$\alpha_{n(q_m)}$ is the right derivative of the concave majorant of
$S^{[0,\infty)}$ at $T_{(n(q_m))} $, Lemma \ref{lem:morethanmu}
implies that $\alpha_{n(q_m)} \te \mu$ almost surely as $m \te
\infty$. This concludes the proof of Theorem \ref{theorem:poisson2}.
\endpf
\vspace{10pt}

\subsection{The structure of the concave majorant of $S^{[0,n]}$ as $n$ varies}
\label{sec:nvaries}

Theorem \ref{theorem:1} relates to the structure of the concave majorant of a random walk of fixed length, and the Theorems \ref{theorem:poisson1} and \ref{theorem:poisson2} allow randomized lengths or infinite length. So far though, we have not discussed how the structure changes as the number of steps of the walk increases, but theorem \ref{theorem:poisson2} and its proof now allow us to make some comments.
Recall that $F_n$ is the number of faces of the concave majorant of $S^{[0,n]} = \{ (j,S_j) \: : \: 0 \le j \le n) \}$, and in the case where $X_1, \ldots , X_n$ are independent with common continuous distribution we know from \re{fnkn} that for each fixed $n$ there is the equality in distribution
$$
F_n \ed K_n := \sum_{j=1}^n I_j 
$$
where the $I_j$ are independent Bernoulli variables with $\P(I_j = 1) = 1/j$.
However, as observed by 
Steele \cite{MR1910531}
the identity in law between $F_n $ and $K_n $ does not hold jointly as $n$ varies, and as
pointed out by 
Qiao and Steele \cite{MR2205729}
the asymptotic behaviour of $F_n $ and $K_n $ as $n \te \infty$ may be quite different.
They provide an example of a continuous distribution of $X_i$ such that for each $m = 1,2, \ldots$
$$
\P( F_n = m \mbox{ infinitely often } ) = 1 
$$
It is an easy consequence of theorem \ref{theorem:poisson2} that
$$
\P( F_n = 1 \mbox{ infinitely often } ) = 1 
$$
if and only if $\E(X^+) = \infty$. It appears that the Poisson analysis of $F_{n(q)}$ can
be used to provide a more thorough description of the possible asymptotic behaviours of
$F_n$ as $n$ varies. In particular, 
as a consequence of the argument of the proof of Lemma \ref{lem:lastvertex}, 
if $\E(X^+) <  \infty$ then $F_n$ is bounded below by the number of faces of the 
majorant on $[0,n]$ which are part of the majorant on $[0,\infty)$, and this number is increasing
in $n$, with limit $\infty$.

\subsection{Decomposition at the maximum}
\label{sec:max}

Theorem \ref{theorem:poisson1}  provides tools for analyzing the behaviour of the random walk $S^{[0,n(q)]}$ before and after the time it achieves its maximum. By conditioning on $n(q)=n$, we can then do the same for $S^{[0,n]}$. The key idea is that by taking the faces of the concave majorant that have positive slope we get only those faces that lie in the region up to where the random walk achieves its maximum, and by taking the faces with negative slope we get only those faces that lie in the region after the time when the random walk achieves its maximum. This approach was used by Spitzer to find identities involving the maximum of a random walk \cite{MR0079851}, as indicated in Section \ref{sec:intro}.


Let $X_1, X_2, \ldots $ be a sequence of independent random variables with common continuous distribution, and let $S_0=0$ and $S_j = \sum_{i=1}^j X_i$ for $j \ge 1 $. Let $S^{[0,n]} = \{ (j,S_j) \: : \: 0 \le j \le n \}$ and $S^{[0,n(q)]} = \{ (j,S_j) \: : \: 0 \le j \le n(q) \}$.
Let $L_n$ be the almost surely unique time at which $S^{[0,n]}$ achieves its maximum, and let the value of the maximum be $M_n$.
Let $F_n$ denote the number of faces of the concave majorant of the walk $S^{[0,n]}$,
with the convention $F_0 = 0$, and let $(N_{n,i}, \Delta_{n,i})$ denote the length and increment associated with the $i$th of these faces.
We make similar definitions when $n$ is randomized to $n(q)$.

\begin{theorem}
\label{theorem:compoundpois}
$(L_{n(q)},M_{n(q)})$ and $(n(q) - L_{n(q)}, S_{n(q)} - M_{n(q)} )$ are independent and both have compound Poisson distributions. 
\end{theorem}

As discussed in Section \ref{sec:intro} the compound Poisson nature of $M_{n(q)}$ and $S_{n(q)}-M_{n(q)}$ and their independence was discovered by Greenwood and Pitman \cite{greenpit80}, but this section gives a more explicit explanation of their distribution.

\vspace{10pt}
\proof
By construction
$$
\Delta_{n,i} = S_{N_{n,1} + \cdots N_{n,i-1} + N_{n,i} } - S_{N_{n,1} + \cdots N_{n,i-1}}
$$
and
$$
\begin{array}{rcl}
(L_{n},M_n) &=& \sum_{i = 1}^{K_n} (N_{n,i} ,\Delta_{n,i})1( \Delta_{n,i} > 0 ) \\
(n - L_{n}, S_n - M_n )  &=& \sum_{i = 1}^{K_n} ( N_{n,i} ,\Delta_{n,i} ) 1( \Delta_{n,i} \le  0 ) 
\end{array}
$$
From Theorem \ref{theorem:poisson1} the $(N_{n(q),i}, \Delta_{n(q),i})$ are the points of a Poisson point process on
$\{ 1,2 \ldots \} \times \mathbb{R}$ with intensity $j^{-1} q^j \P( S_j \in dx) , j \in \{ 1,2,\ldots \}, x \in \mathbb{R}$, and thus the conclusion follows.
\endpf
\vspace{10pt}


In the special case where $\P(S_j > 0)$ is constant for $1 \le j \le
n$, by conditioning on the event $n(q) = n$ and $L_{n(q)} = \ell$ we
can deduce results about the concave majorant of $S^{[0,n]}$ either
side of its maximum.
\begin{theorem}\label{theorem:maxdecomp}
Let $X_1,\ldots,X_n$ be independent with common continuous distribution. Let $S_0=0$ and $S_j = \sum_{i=1}^j X_i$ for $1 \le j \le n $, and let $S^{[0,n]} = \{ (j,S_j) \: : \: 0 \le j \le n \}$. Suppose that $\P(S_j > 0) = p_+$ for $1 \le j \le n$. Then conditionally given $L_n : = \arg \max_{0 \le j \le n}S_j = \ell$, the partition generated by the lengths of the faces of the concave majorant of $S^{[0,n]}$ on the interval $[0,\ell]$ is distributed according to the Ewens sampling formula with parameter $p_+$.
That is, if $A_j^+$ is the number of faces of the concave majorant with positive slope of length $j$, then for any $\{ a_j \: : \: j \geq 1 \}$ such that $\sum_j j a_j = \ell \leq n$,
\eq
\lb{eq:ewens}
\P(A_j^+ = a_j , j \ge 1 | L_n = \ell) = \frac{ \Gamma(p_+) \ell! }{ \Gamma(p_+ + \ell) } \prod_{j=1}^{\ell} \frac{ (p_+)^{a_j} }{ j^{a_j} a_j ! }
\en
The partition generated by the lengths of the faces of the concave majorant of $S^{[0,n]}$ on the interval $[\ell,n]$ is also distributed according to the Ewens sampling formula but with parameter $p_-=1-p_+$. 
\end{theorem}
\proof
Let $A_{n(q),j}^+$ be the number of faces of the concave majorant of $S^{[0,n(q)]}$ with positive slope of length $j$. From the proof of Theorem \ref{theorem:compoundpois} it is easy to see that $A_{n(q),j}^+$ has a Poisson distribution with parameter $j^{-1}q^j p_-$, independently for each $j$, and independently of $S^{[0,n(q)]}$ after time $L_{n(q)}$. Thus for any $\{ a_j \: : \: j \geq 1 \}$ such that $\sum_j j a_j = \ell$,
\begin{eqnarray}
\P(A_j^+ = a_j , j \ge 1 | L_n = \ell) 
& = & \P(A_{n(q),j}^+ = a_j , j \ge 1 | L_{n(q)}= \ell , n(q) = n ) \nonumber \\
& = & \P(A_{n(q),j}^+ = a_j , j \ge 1 | L_{n(q)} = \ell ) \nonumber  \\
& = & \frac{  \P ( A_{n(q),j}^+ = a_j , j \ge 1 ) }{\P( L_{n(q)} = \ell )} \nonumber  \\
& = & \frac{ \prod_j \frac{ (p_+)^{a_j} q^{ja_j}} { j^{a_j} a_j ! } \exp \{ -\frac{ p_+ q^j }{ j } \} } {\P( L_{n(q)} = \ell )} \label{eq:preewens}
\end{eqnarray}
Under the assumption $\P(S_j > 0) = p_+$ for $1 \le j \le n$, it is known \cite[Chapter XII, ($8.12$)]{MR0270403} that for the random walk $S^{[0,n]}$, the almost surely unique index $L_n$ such that $S_{L_n} = \max_{0 \le j \le n}S_j$ has the beta-binomial distribution
$$
\P( L_n = \ell ) = (-1)^n { p_- - 1 \choose \ell } { p_+ -1 \choose n - \ell }   ~~~~ (0 \le \ell \le n )
$$
which is the mixture of binomial$(n,p)$ distributions for $p$ with beta$(p_+,p_-)$ distribution on $[0,1]$. Thus
$$
\P( L_{n(q)} = \ell ) = \frac{ \Gamma( p_+ + \ell ) q^{\ell} (1-q)^{p_+} }{ \Gamma( p_+ ) \ell !  } 
$$
Thus \re{eq:preewens} reduces to \re{eq:ewens}. The partition after the maximum is proved similarly.
\endpf

\section{The general case}
\label{sec:noncontinuous}

Let $S_j = \sum_{i=1}^j X_i$ for $1 \le j \le n$, where $X_1, X_2, \ldots $ is a sequence of exchangeable random variables. Let $S^{[0,n]} = \{ (j,S_j) \: : \: 1 \le j \le n \}$, and let $\bar{C}^{[0,n]}$ be the concave majorant of $S^{[0,n]}$. The concave majorant in this case, where there may some subsets of $X_1,\ldots,X_n$ that have the same arithmetic mean, is less well studied. However, the literature does contain some results for the case where $X_1,X_2,\ldots $ are also assumed to be independent.

Sparre Andersen \cite{andersenH_n} introduced the random variable $H_n$, the number of $1 \le j \le n$ such that $S_j = \bar{C}^{[0,n]}(j)$, and $F_n$, the number of faces of the concave majorant, i.e.\ the number of distinct slopes in the concave majorant (note that Andersen uses $K_n$ instead of $F_n$, but we will always use $K_n$ to represent the number of cycles in a random permutation of $[n]$). Figure \ref{fig:walk} shows an example of a random walk with $F_n=3$ and $H_n=8$. Clearly, $F_n \le H_n$, and in the case of continuous distributions we have $F_n = H_n$ almost surely. Sparre Andersen derived the generating function
\eq
\lb{Hgen}
H(s,t) := \sum_{n=0}^{\infty} \sum_{m=0}^{n} \P(H_n = m) s^n t^m 
\en
for all distributions of $X_1$.
As will be shown in Theorem \ref{theorem:HKFgen} the theory presented in this section provides a powerful new method of deriving this formula, and in addition a formula for a similar generating function involving $F_n$.

Sherman \cite{MR0169306} introduced a further variable $J_n$ relating to the concave majorant with $H_n \leq J_n \leq F_n$. Sherman deduces a Spitzer identity which relates the generating functions of $J_n$ and $\Phi_n$, the periodicity of $(X_1, \ldots , X_n)$, that is, the maximal number $\phi$ such that $(X_1, \ldots , X_n) = (X_1, \ldots , X_{n / \phi} , \: \: \ldots \: \: , X_1 , \ldots , X_{n / \phi})$.

In this section it will be important to make a distinction between \emph{excursions}, \emph{segments} and \emph{faces}, and between their associated compositions of $n$. The following definitions are illustrated in Figure \ref{fig:walk}.
\begin{itemize}
\item
An \emph{excursion} is a section of a walk between two integer valued times with the property that the walk touches its concave majorant at the end points of the excursion but lies strictly below it between the end points. The number of distinct excursions of $S^{[0,n]}$ is equal to $H_n$. Let $\Xi_{[0,n]}^H$ be the composition of $n$ induced by the lengths of the excursions of $S_{\rho}^{[0,n]}$, the transformed walk of Theorem \ref{theorem:1}. Although this has the same distribution as the composition induced by the lengths of the excursions of $S^{[0,n]}$, the forthcoming discussion about \emph{segment} compositions only makes sense for $S_{\rho}^{[0,n]}$. We say that the \emph{slope} of an excursion is the slope of the line joining its start and end points.
\item
A \emph{segment} will always refer to one segment of a partition. That is, if $(n_1,\ldots, n_k)$ a partition of $n$ then we say it has $k$ segments with associated lengths $n_1,\ldots,n_k$. 
As we described in the introduction, to generate a walk with the law of $S^{[0,n]}$ whilst simultaneously getting information about its concave majorant, i.e.\ to generate $S_{\rho}^{[0,n]}$, we first choose a random partition induced by the cycle lengths of a uniform random permutation. If we are just interested in the concave majorant of $S_{\rho}^{[0,n]}$, then we only need to associate a slope with each segment of that partition and then arrange the segments in order of non-increasing slope, where the ordering of any segments with the same slope is chosen uniformly randomly. Keeping track of the end points of the segments results in another induced composition of $n$, which we call $\Xi_{[0,n]}^K$. This composition arises from our construction and cannot be read off from a given random walk.
\item
A \emph{face} will mean one face of the concave majorant. The number of distinct faces is equal to $F_n$. Let $\Xi_{[0,n]}^F$ be the composition of $n$ induced by the lengths of the faces of $S_{\rho}^{[0,n]}$. Again, this has the same distribution as the composition of $n$ induced by the lengths of the faces of $S^{[0,n]}$.
\item
The terms \emph{excursion block}, \emph{segment block} and \emph{face block} will mean blocks of the compositions $\Xi_{[0,n]}^H$, $\Xi_{[0,n]}^K$ and $\Xi_{[0,n]}^F$ respectively, where for example the \emph{blocks} of the composition $(3,4,1)$ of $8$ in order are defined to be $[0,3]$, $[3,7]$ and $[7,8]$. The slope associated with any block $[a,b]$ is defined by $(S^{\rho}_b - S^{\rho}_a)/(b-a)$.
\end{itemize}

\begin{figure}
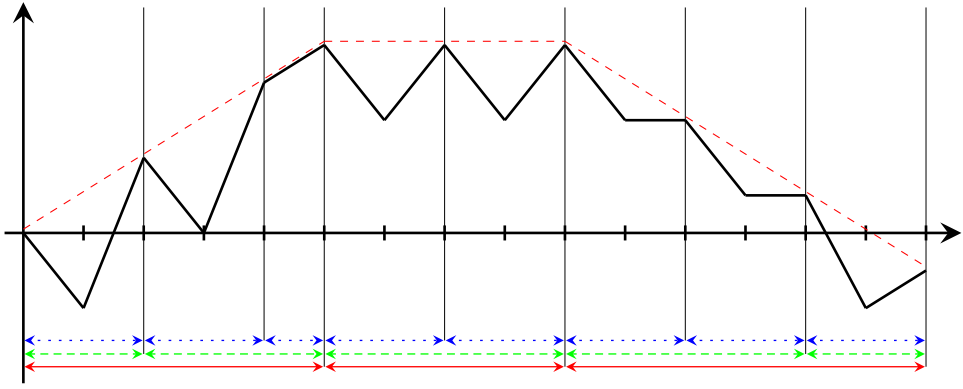

\begin{pgfpicture}{0.7cm}{0cm}{19cm}{5.02cm}
\begin{pgfscope} 
\pgfsetendarrow{\pgfarrowsingle}
\pgfsetlinewidth{1pt} 
\pgfxyline(0.55,2)(13.2,2) 
\pgfxyline(0.8,0)(0.8,5) 
\end{pgfscope} 
\begin{pgfscope} 
\pgfsetlinewidth{1pt} 
\pgfxyline(1.6,1.9)(1.6,2.1) 
\pgfxyline(2.4,1.9)(2.4,2.1) 
\pgfxyline(3.2,1.9)(3.2,2.1) 
\pgfxyline(4.0,1.9)(4.0,2.1) 
\pgfxyline(4.8,1.9)(4.8,2.1) 
\pgfxyline(5.6,1.9)(5.6,2.1) 
\pgfxyline(6.4,1.9)(6.4,2.1) 
\pgfxyline(7.2,1.9)(7.2,2.1) 
\pgfxyline(8.0,1.9)(8.0,2.1) 
\pgfxyline(8.8,1.9)(8.8,2.1) 
\pgfxyline(9.6,1.9)(9.6,2.1) 
\pgfxyline(10.4,1.9)(10.4,2.1) 
\pgfxyline(11.2,1.9)(11.2,2.1) 
\pgfxyline(12.0,1.9)(12.0,2.1) 
\pgfxyline(12.8,1.9)(12.8,2.1) 
\end{pgfscope} 
\begin{pgfscope}
\pgfsetlinewidth{0.1pt}
\pgfxyline(2.4,0.39)(2.4,5)
\pgfxyline(4.0,0.57)(4.0,5)
\pgfxyline(4.8,0.22)(4.8,5) 
\pgfxyline(6.4,0.57)(6.4,5)
\pgfxyline(8.0,0.22)(8.0,5)
\pgfxyline(9.6,0.57)(9.6,5)
\pgfxyline(11.2,0.39)(11.2,5)
\pgfxyline(12.8,0.22)(12.8,5)
\end{pgfscope}
\begin{pgfscope}
\color{blue}
\pgfsetlinewidth{0.5pt}
\pgfsetdash{{1pt}{3pt}}{0pt}
\pgfsetendarrow{\pgfarrowto}
\pgfsetstartarrow{\pgfarrowto}
\pgfsetendarrow{\pgfarrowsingle}
\pgfsetstartarrow{\pgfarrowsingle}
\pgfxyline(0.85,0.57)(2.35,0.57)
\pgfxyline(2.45,0.57)(3.95,0.57) 
\pgfxyline(4.05,0.57)(4.75,0.57) 
\pgfxyline(4.85,0.57)(6.35,0.57)
\pgfxyline(6.45,0.57)(7.95,0.57) 
\pgfxyline(8.05,0.57)(9.55,0.57)
\pgfxyline(9.65,0.57)(11.15,0.57)
\pgfxyline(11.25,0.57)(12.75,0.57) 
\end{pgfscope}
\begin{pgfscope}
\color{green}
\pgfsetlinewidth{0.5pt}
\pgfsetdash{{3pt}{2pt}}{0pt}
\pgfsetendarrow{\pgfarrowto}
\pgfsetstartarrow{\pgfarrowto}
\pgfsetendarrow{\pgfarrowsingle}
\pgfsetstartarrow{\pgfarrowsingle}
\pgfxyline(0.85,0.39)(2.35,0.39)
\pgfxyline(2.45,0.39)(4.75,0.39) 
\pgfxyline(4.85,0.39)(7.95,0.39)
\pgfxyline(8.05,0.39)(11.15,0.39)
\pgfxyline(11.25,0.39)(12.75,0.39) 
\end{pgfscope}
\begin{pgfscope}
\color{red}
\pgfsetlinewidth{0.5pt}
\pgfsetendarrow{\pgfarrowto}
\pgfsetstartarrow{\pgfarrowto}
\pgfsetendarrow{\pgfarrowsingle}
\pgfsetstartarrow{\pgfarrowsingle}
\pgfxyline(0.85,0.22)(4.75,0.22) 
\pgfxyline(4.85,0.22)(7.95,0.22)
\pgfxyline(8.05,0.22)(12.75,0.22)
\end{pgfscope}
\begin{pgfscope}
\color{red}
\pgfsetdash{{3pt}{3pt}}{0pt}
\pgfxyline(0.8,2.05)(4.8,4.55)
\pgfxyline(4.8,4.55)(8,4.55) 
\pgfxyline(8,4.55)(12.8,1.55) 
\color{black}
\end{pgfscope}
\begin{pgfscope} 
\pgfsetlinewidth{1pt} 
\pgfxyline(0.8,2)(1.6,1) 
\pgfxyline(1.6,1)(2.4,3) 
\pgfxyline(2.4,3)(3.2,2) 
\pgfxyline(3.2,2)(4.0,4) 
\pgfxyline(4.0,4)(4.8,4.5) 
\pgfxyline(4.8,4.5)(5.6,3.5) 
\pgfxyline(5.6,3.5)(6.4,4.5) 
\pgfxyline(6.4,4.5)(7.2,3.5) 
\pgfxyline(7.2,3.5)(8.0,4.5) 
\pgfxyline(8.0,4.5)(8.8,3.5) 
\pgfxyline(8.8,3.5)(9.6,3.5) 
\pgfxyline(9.6,3.5)(10.4,2.5) 
\pgfxyline(10.4,2.5)(11.2,2.5) 
\pgfxyline(11.2,2.5)(12.0,1.0) 
\pgfxyline(12.0,1.0)(12.8,1.5) 
\end{pgfscope} 
\end{pgfpicture}
\captionsetup{singlelinecheck=off}
\caption{An example of a random walk with non-continuous increment
  distribution for $n=15$, with $F_n = 3$ and $H_n = 8$. The concave
  majorant is shown with dashed line. The compositions induced by the excursion
  lengths and face lengths are fixed by the values of the walk, and an
  example of a possible composition induced by the lengths of the chords associated with the partition segments is shown. The
  compositions going from top to bottom are $\Xi^H_{[0,15]} = (2,2,1,2,2,2,2,2,)$, $\Xi^K_{[0,15]} =
  (2,3,4,4,2)$ and $\Xi^F_{[0,15]} = (5,4,6)$.}
\label{fig:walk}
\end{figure}

Since the values of any walk on $[0,n]$ between two vertices of its concave majorant, i.e.\ between the start and end points of some face, are composed of one or many consecutive excursions, $ \Xi_{[n]}^H$ is some refinement of $\Xi_{[n]}^F$, which we write as $ \Xi_{[n]}^H \preceq \Xi_{[n]}^F$.
For $S_{\rho}^{[0,n]}$ constructed as in Theorem \ref{theorem:1}, define $H_n^{\rho}$ and $F_n^{\rho}$ similarly to $H_n$ and $F_n$, and note that $H_n \ed H_n^{\rho}$ and $F_n \ed F_n^{\rho}$. Recall that $K_n$ is the number of segments in the partition chosen at the beginning of the construction. We will have $H_n^{\rho} \le K_n \le F_n^{\rho}$, and moreover $\Xi_{[0,n]}^K$ will be such that $\Xi_{[0,n]}^H \preceq \Xi_{[0,n]}^K \preceq \Xi_{[0,n]}^F$. We will discuss these nested compositions further after proving Theorem \ref{theorem:1} in the general case.

\vspace{10pt}
\noindent{\textbf{Proof. (Theorem \ref{theorem:1})} As in the proof of Theorem \ref{theorem:1} under assumption $\mathbf{A}$, it is enough to show that if $X_1,\ldots,X_n$ are samples without replacement from a list $x_1,\ldots,x_n$ of real numbers, where now each number is labelled but no longer necessarily distinct in value, then
$$
\P(X_{\rho(1)} = x_1, \ldots , X_{\rho(n)} = x_n ) = \frac{1}{n!}
$$
Let $x=(x_1,\ldots,x_n)$, and suppose this is fixed throughout the proof of the theorem. Let $\bar{c}^{[0,n]}$ be the concave majorant of the deterministic walk with increments $x_1,\ldots,x_n$. Some notation and a couple of combinatorial lemmas are needed before continuing. 

For any $n \in \mathbb{N}$, let $\mathcal{N}_n$ be the set of all compositions of $n$. Let $f \in \mathbb{N}$, $h \in \mathbb{N}$ and $(v_1,\ldots,v_f) \in \mathcal{N}_h$. Let  $\mathcal{N}_{(v_1,\ldots,v_f) , (k_1,\ldots,k_f)}$ be the set
$$
\{ (h_1,\ldots,h_{\sum_{i=1}^f k_i}) \in \mathcal{N}_h  \: : \:   (h_{\sum_{i=1}^{j-1} k_i}, \ldots , h_{\sum_{i=1}^j k_i}) \in \mathcal{N}_{v_j} \text{ for } 1 \le j \le f \}
$$
Thus an element of $\mathcal{N}_{(v_1,\ldots,v_f) , (k_1,\ldots,k_f)}$ is a composition of $h$ formed by joining together compositions of $v_1,\ldots,v_f$ which contain $k_1,\ldots,k_f$ blocks respectively (and hence $\mathcal{N}_{(v_1,\ldots,v_f) , (k_1,\ldots,k_f)} $ may be an empty set for some values of $(k_1,\ldots,k_f)$). 

\begin{lemma}\label{lem:stirlingnumbers}
Let $f \in \mathbb{N}$, $h \in \mathbb{N}$ and $(v_1,\ldots,v_f) \in \mathcal{N}_h$. Then
\eq
\lb{bigsum}
\sum_{k=f}^h \sum_{(k_1,\ldots,k_f) \in \mathcal{N}_k } \sum_{ (h_1 , \ldots ,h_k) \in \mathcal{N}_{(v_1,\ldots,v_f),(k_1,\ldots,k_f)}  } \prod_{i=1}^k \frac{1}{k_1! \cdots k_f!} \frac{1}{h_1 \cdots h_k} = 1
\en
\end{lemma}
\proof
The numbers that are being summed over bear a strong resemblance to the unsigned Stirling numbers of the first kind $|S(n,k)|$, which enumerate the number of permutations of $n$ with $k$ cycles. Using this as a guide, consider a set $A$ consisting of permutations of $v_1,\ldots,v_f$, where permutations corresponding to $v_i$ and $v_j$ with $i \neq j$ are considered distinct even if they are identical. The number of such sets where for each $1 \le j \le f$ the permutation of $v_j$ has $k_j$ cycles of sizes $h_{\sum_{i=1}^{j-1} k_i}, \ldots , h_{\sum_{i=1}^j k_i}$ is
$$
\frac{v_1! \cdots v_f!}{k_1! \cdots k_f! \: \cdot \: h_1 \cdots h_k } 
$$
Since the total number of elements of $A$ is $v_1! \cdots v_f!$, and the summation in \re{bigsum} simplifies to be the sum over the subsets of $A$ such that for each $1 \le j \le f$ the permutation of $v_j$ has $k_j$ cycles of sizes $h_{\sum_{i=1}^{j-1} k_i}, \ldots , h_{\sum_{i=1}^j k_i}$, the value of the sum must be 1.
\endpf
\vspace{10pt}

Let $f(\bar{c}^{[0,n]})$ be the number of faces of $\bar{c}^{[0,n]}$, and let
$\ell_1(\bar{c}^{[0,n]}), \ldots , \ell_{f(\bar{c}^{[0,n]} )}(\bar{c}^{[0,n]})$
be the lengths of those faces, arranged in the order those faces appear in $\bar{c}^{[0,n]}$.
Let $\mathcal{N}(\bar{c}^{[0,n]})$ be the set
$$
\! \! \! \! \! \! \! \! \{ (n_1,\ldots,n_k) \in \mathcal{N}_n \: :  \: \:  \exists \mbox{ } k_1 < \cdots < k_{f(\bar{c}^{[0,n]})} \text{ s.t. } \sum_{i=k_{j-1}}^{k_j} n_i = \ell_j(\bar{c}^{[0,n]}), 1 \le j \le f(\bar{c}^{[0,n]}) \}  
$$
Loosely, $\mathcal{N}(\bar{c}^{[0,n]}) $ is the set of possible values for $\Xi_{[0,n]}^K$ conditionally given that the concave majorant of $S_{\rho}^{[0,n]}$ is $\bar{c}^{[0,n]}$. For $(n_1,\ldots, n_k) \in \mathcal{N}(\bar{c}^{[0,n]})$, let
$$
\{ k_j(n_1,\ldots,n_k) , 1 \le j \le f(\bar{c}^{[0,n]}) \} 
= 
\{ (k_1,\ldots,k_{f(\bar{c}^{[0,n]})} )  \: : \: 
\sum_{i=k_{j-1}}^{k_j} n_i = \ell_j(\bar{c}^{[0,n]}) \}
$$
Then $k_j(\Xi_{[0,n]}^K)$ represents the number of blocks of $\Xi_{[0,n]}$ that lie in the $j$th face block, i.e.\ in the $j$th block of $\Xi_{[0,n]}^F$.
Finally, let
$$
\mathcal{N}_x(\bar{c}^{[0,n]}) = \{ (n_1,\ldots,n_k) \in \mathcal{N}(\bar{c}^{[0,n]})   \: : \: \: \sum_{j=1}^{n_i} x_j = \bar{c}^{[0,n]}(n_i) \text{ for } 1 \le i \le k \}  
$$
Then $\mathcal{N}_x(\bar{c}^{[0,n]}) $ is the set of possible values for $\Xi_{[0,n]}^K$ conditionally given that 
$\{ X_{\rho(i)} = x_i \: : \: 1 \le i \le n \}$.

\begin{lemma}\label{lem:sumsto1}
For every composition $(n_1,\ldots,n_k) \in \mathcal{N}_x(\bar{c}^{[0,n]})$, for $1 \le i \le k$ let
$$
h_i(x,n_1,\ldots,n_k) = \# \{ j : n_1 + \cdots + n_{i-1} < j \le n_1 + \cdots + n_i \, , \, \mbox{$\sum_{l=1}^j x_l$} = \bar{c}^{[0,n]}(j) \}
$$
Then
\eq
\lb{compsums}
\sum_{k=1}^n \sum_{(n_1,\ldots,n_k) \in \mathcal{N}_x(\bar{c}^{[0,n]}) }
\left( \prod_{i=1}^k \frac{1}{ h_i(x,n_1,\ldots,n_k) } \right)  
\left( \prod_{j=1}^{f(\bar{c}^{[0,n]} )} \frac{1}{k_j(n_1,\ldots,n_k)!} \right) 
= 1
\en
\end{lemma}
\proof
Let $h = \# \{ j : 1 \le j \le n, \sum_{l=1}^j = \bar{c}^{[0,n]}(j) \}$ and for $1 \le i \le f(\bar{c}^{[0,n]})$ let
\begin{eqnarray*}
&& \! \! \! \! \! \!  \! \! \!  \! \! \!  v_i(x)   =   \\
&& \! \! \! \! \! \! \! \! \! \! \! \! \! \! \! \! \! \! \! \! \# \{ j : \ell_1(\bar{c}^{[0,n]}) + \cdots + \ell_{i-1}(\bar{c}^{[0,n]}) < j \le \ell_1(\bar{c}^{[0,n]}) + \cdots + \ell_i(\bar{c}^{[0,n]}) \, , \, \mbox{$\sum_{l=1}^j x_l$} = \bar{c}^{[0,n]}(j) \}  
\end{eqnarray*}
Associate with each composition $(n_1, \ldots , n_k) \in \mathcal{N}_x(\bar{c}^{[0,n]})$ of length $k$
a composition of $h$
$$
(h_1(x,n_1,\ldots,n_k) , h_2(x,n_1,\ldots,n_k) , \ldots , h_k(x,n_1,\ldots,n_k) )
$$
so that there is a bijection between the elements of $\mathcal{N}_x(\bar{c}^{[0,n]})$ with $k$ blocks and the set of compositions $(h_1,\ldots,h_k)$ of $h$ with $k$ blocks that are formed by joining together in order compositions of $v_1,\ldots,v_{f(\bar{c}^{[0,n]})}$ which have $k_1,\ldots,k_{f(\bar{c}^{[0,n]})}$ blocks respectively. Thus the term on the left hand side of \re{compsums} is 
$$
\sum_{k=f}^h \sum_{(k_1,\ldots,k_{f(\bar{c}^{[0,n]})}) \in \mathcal{N}_k } \sum_{ (h_1 , \ldots h_k) \in \mathcal{N}_{(v_1,\ldots,v_f),(k_1,\ldots,k_f)}  } \prod_{i=1}^k \frac{1}{k_1! \cdots k_f!} \frac{1}{h_1 \cdots h_k} 
$$
which by Lemma \ref{lem:stirlingnumbers} is 1.
\endpf
\vspace{10pt}

Fix a composition $(n_1,\ldots,n_k)$ of $n$. For $1 \le j \le n$ let $I_j = \{ i \: : \: n_i = j \}$ and let $a_j = |I_j|$. Following the construction of $S_{\rho}^{[0,n]}$ described in the introduction, we see that the event $\{ \Xi_{[0,n]}^K = (n_1,\ldots,n_k) \text{ and } X_{\rho(\ell)} = x_{\ell}, 1 \le \ell \le n \}$ occurs if and only if
\begin{enumerate}[(i)]
\item
$L_{n,1},\ldots,L_{n,K_n}$ is $(n_1,\ldots,n_k)$ in non-increasing order;
\item 
for each $1 \le j \le n$, for each $i \in I_j$ the ordered list $(X_{n_1+\cdots+n_{i-1}+1}, \ldots, X_{n_1+\cdots+n_i})$ is one of the $n_i=j$ cyclic permutations of the ordered list 
\newline $( x_{m_1+m_2+\cdots+m_{\tau(i')-1}+1}, \ldots , x_{m_1+m_2+\cdots+m_{\tau(i')}} ) $ for some $i' \in I_j$;
\item
for each $1 \le j \le n$, for each $i \in I_j$ the cyclic permutation that is chosen for the ordered list of increments $(X_{n_1+\cdots+n_{i-1}+1}, \ldots, X_{n_1+\cdots+n_i})$ is the unique cyclic permutation that results in the ordered list becoming exactly $( x_{m_1+m_2+\cdots+m_{\tau(i')-1}+1}, \ldots , x_{m_1+m_2+\cdots+m_{\tau(i')}} ) $;
 \item
for each $1 \le j \le f(\bar{c})^{[0,n]})$ the ordering of the $k_j(n_1,\ldots,n_k)$ segments within the $j$th face is chosen correctly out of the $k_j!$ possible orderings.
\end{enumerate} 
Recall that  for $1 \le i \le k$ we have
$$
h_i(x,n_1,\ldots,n_k) = \# \{ j : n_1 + \cdots + n_{i-1} < j \le n_1 + \cdots + n_i \, , \, \mbox{$\sum_{l=1}^j x_l$} = \bar{c}^{[0,n]}(j) \}
$$
so that in (iii) there are $\prod_{i=1}^k h_i(x,n_1,\ldots,n_k)$ possible choices of combinations of cyclic permutations. 
Then the probability of the event 
\newline $\{ \Xi_{[0,n]}^K = (n_1,\ldots,n_k) \text{ and } X_{\rho(\ell)} = x_{\ell}, 1 \le \ell \le n \}$ is
$$
\! \! \! \! \! \! \! \! \! \! 
\left(\prod_{j=1}^n \frac{1}{a_j!} \: \prod_{i=1}^k \frac{1}{n_i} \right) \left( \frac{1}{n!} \prod_{i=1}^k n_i \: \prod_{j=1}^n a_j! \right) 
\left( \prod_{i=1}^k \frac{1}{ h_i(x,n_1,\ldots,n_k) } \right)  
\left( \prod_{j=1}^{f(\bar{c}^{[0,n]} )} \frac{1}{k_j(n_1,\ldots,n_k)!} \right) 
$$
where the first two terms should be familiar from the proof of Theorem \ref{theorem:1} under assumption $\mathbf{A}$. Finally, by summing this probability over all possible compositions, we have that the probability of the event
$ \{ X_{\rho(\ell)} = x_{\ell}, 1 \le \ell \le n \}$ is 
$$
\frac{1}{n!} \sum_{k=1}^n \sum_{(n_1,\ldots,n_k) \in \mathcal{N}_x(\bar{c}^{[0,n]}) }
\left( \prod_{i=1}^k \frac{1}{ h_i(x,n_1,\ldots,n_k) } \right)  
\left( \prod_{j=1}^{f(\bar{c}^{[0,n]} )} \frac{1}{k_j(n_1,\ldots,n_k)!} \right) 
= \frac{1}{n!}
$$
where the equality is by Lemma \ref{lem:sumsto1}. 
This completes the proof of Theorem \ref{theorem:1}.
\endpf
\vspace{10pt}

In the case where $X_1,X_2,\ldots$ are independent, the Poisson point process ideas of Section \ref{sec:poisson} lead to a simpler description of the concave majorant. For the rest of this section it is assumed that $X_1, X_2, \ldots $ is a sequence of  independent and identically distributed random variables and $n(q)$ is a geometric variable with parameter $1-q$. Let $S^{[0,n(q)]} = \{ (j,S_j) : 0 \le \j \le n(q) \}$, where $S_0=0$ and  $S_j = \sum_{i=1}^j X_i$ for $j \ge 1$. Let $\bar{C}^{[0,n]}$ be the concave majorant of $S^{[0,n(q)]}$. The following theorem is the extension to the non-continuous increment case of Theorem \ref{theorem:poisson1}.

\begin{theorem}\label{theorem:poissonnoncont}
If $X_1, X_2, \ldots $ are independent with common distribution and $n(q)$ a geometric variable with parameter $1-q$, then the lengths and increments of the faces of the concave majorant of the random walk $S^{[0,n(q)]}$ have the following law. Let $\mathfrak{P}$ be a Poisson point process of on $\{ 1,2,\ldots \} \times \mathbb{R}$ with intensity $j^{-1} q^j \P(S_j \in dx)$ for $j = 1, 2, \ldots$, $x \in \mathbb{R}$. Note that this process may result in multiple points at the same location. Each point of $\mathfrak{P}$ represents the length and increment of a chord associated with some segment of a partition of $n(q)$. Chords with the same slope are joined together in uniform random order, independently of their lengths, to form the faces of the concave majorant. Moreover, let $K_{n(q)}$ be the total number of chords associated with partition segments and for $1 \le i \le K_{n(q)}$ let $N_{n(q),i}$ be the length of the $i$th of these chords once they have been ordered by decreasing slope and uniform randomization of ties. Then the sequence of path segments
$$
\{ (S_{\sum_{l = 1}^{i-1} N_{n(q),l} + k} - S_{\sum_{l = 1}^{i-1} N_{n(q),l} },  0 \le k \le N_{n(q),i} ) , i = 1, \ldots , K_{n(q)} \}
$$
is a list of the points of a Poisson point process in the space of finite random walk segments
$$
\{ (s_1, \ldots , s_j) \mbox{ for some } j = 1,2, \ldots \}
$$
whose intensity measure on paths of length $j$ is $j^{-1}$ times the conditional distribution of $(S_1, \ldots , S_j)$ given that $S_k < (k/j) S_j $ for all $1 \le k < j $. Again, this Poisson point process may result in multiple points at the same location.
\end{theorem}
\proof
For any $n \in \mathbb{N}$, conditionally given $n(q)=n$, the projection of the points of $\mathfrak{P}$ onto $\{ 1, \ldots , n \}$ has the law of a partition of $n$ generated by the cycle lengths of a random permutation of $[n]$ by Lemma \ref{lem:poisson1}. Hence we know from Theorem \ref{theorem:1} that for every $n \in \mathbb{N}$, conditionally given $n(q) = n$, the process described in the theorem gives the correct law for the concave majorant of $S^{[0,n]}$ and gives the correct law for $\Xi_{[0,n]}^K$, the composition induced by the lengths of the partition segments involved in creating $S_{\rho}^{[0,n]}$. The remaining assertions follow by independence of the walks associated with each partition segment.
\endpf
\vspace{10pt}

We now move towards describing the joint law of the nested compositions $\Xi_{[0,n(q)]}^H \preceq \Xi_{[0,n(q)]}^K \preceq \Xi_{[0,n(q)]}^F$ in the case where $X_1,X_2,\ldots$ are independent and the walk has geometric length. The full description of this law will be given in Theorem \ref{theorem:HKFcomposition} at the end of this section, along with some applications of the theory.
Let $S_{\rho}^{[0,n(q)]}$ be such that  conditionally given $n(q) = n$, $S_{\rho}^{[0,n(q)]}$ is constructed in the same way as $S_{\rho}^{[0,n]}$ in Theorem \ref{theorem:1}, and let $\bar{C}_{\rho}^{[0,n(q)]}$ be the concave majorant of $S_{\rho}^{[0,n]}$. We begin by describing the laws of $H_{n(q)}$, $K_{n(q)}$ and $F_{n(q)}$, which are defined to be the number of excursions, segments and faces  respectively of $\bar{C}_{\rho}^{[0,n(q)]}$.

We need some new notation, some of which is taken from Sparre Andersen \cite{andersenH_n}. Let $x_1,x_2,\ldots$ be an enumeration of the set of real numbers $x$ for which $\P(S_k=kx)$ is positive for some $k>0$, and let
\begin{eqnarray*}
\mu_j(q) & = & \sum_{k=1}^{\infty} k^{-1} q^k \P(S_k=kx_j), \quad \text{ for } j=1,2,\ldots \\
& & \\
\mu_0(q) & = & \sum_{k=1}^{\infty} k^{-1} q^k \P(S_k \neq kx_j \text{ for } j=1,2,\ldots ) \\
& = & - \log (1-q) - \sum_{j=1}^{\infty} \mu_j(q)
\end{eqnarray*}

\begin{proposition}\label{prop:HKFdist}
Let $H_{q,j}$, $K_{q,j}$ and $F_{q,j} $ be the number of excursion, segments and faces in $\bar{C}_{\rho}^{[0,n(q)]}$ of slope $x_j$ for $j \ge 1$. Then for each $j \ge 1$
\begin{enumerate}[(i)]
\item
$H_{q,j}$ is a geometric random variable with parameter $\exp ( -\mu_j(q) )$, independently of $\{ H_{q,i}:i \neq j \}$.
\item
$K_{q,j}$ is a Poisson random variable with parameter $\mu_j(q)$, independently of $\{ K_{q,i} :i \neq j \}$.
\item
$F_{q,j}$ is a Bernoulli random variable with parameter $1 - \exp ( -\mu_j(q) )$, independently of $\{ F_{q,i}:i \neq j \}$.
\end{enumerate} 
Let $H_{q,0}$, $K_{q,0}$ and $F_{q,0} $ be the number of excursion, segments and faces with slope not equal to $x_j$ for any $j \ge 1$. Then
\begin{enumerate}[(i)]
  \setcounter{enumi}{3}
  \item $H_{q,0}=K_{q,0}=F_{q,0}$ almost surely and their common distribution is Poisson with parameter $\mu_0(q)$, independently of $\{ H_{q,j},K_{q,j},F_{q,j}: j \ge 1 \}$.
\end{enumerate}
\end{proposition}
\proof
(ii) follows from Theorem \ref{theorem:poissonnoncont}, (iii) is implied by (ii) since a face of slope $x$ exists if and only if there is at least one segment of slope $x$, and (iv) is also implied by Theorem \ref{theorem:poissonnoncont} since it concerns the restriction of the Poisson point process to slopes which have zero probability, as in the case of continuous increment distributions.

Fix $j \ge 1$. (ii) implies that $\P(H_{q,j} \geq 1) = \P(K_{q,j} \geq 1) = 1-\exp ( -\mu_j(q) )$. Given that there at least $n$ excursions of slope $x_j$, by the memoryless property of the geometric distribution of $n(q)$, the law of the remaining values of the walk $S_{\rho}^{[0,n(q)]}$ is the same as the law of a walk generated by the Poisson process of path segments in Theorem \ref{theorem:poissonnoncont} but thinned to only include segments with slope $x \ge x_j$.
Thus
$$
\P( H_{q,j} \geq n+1 | H_{q,j} \geq n ) = \P(K_{q,j} \ge 1) = 1-\exp ( -\mu_j(q) )
$$
which proves (i).
\endpf

\begin{theorem}\label{theorem:HKFgen}
Let $H_n$ and $F_n$ be the number of excursions and faces for $S^{[0,n]}$, and let $K_n$ be the number of segments for $S_{\rho}^{[0,n]}$. 
Then for $0 \le s,t \le 1$,
\begin{eqnarray*}
H(s,t) & = & e^{ t \mu_0(s) } \prod_{j=1}^{\infty} \frac{1}{ 1-t + te^{-\mu_j(s)}} \\
K(s,t) & = & e^{ t \mu_0(s) } \prod_{j=1}^{\infty} e^{t \mu_j(s)} \: \: =  \: \: (1-s)^{-t}  \\
F(s,t) & = &  e^{ t \mu_0(s) } \prod_{j=1}^{\infty} (1-t + te^{\mu_j(s)})
\end{eqnarray*}
\end{theorem}

The generating function of $G_{K_n}(z) = \sum_{m=1}^{\infty} z^m \P(K_n=m) $ is well known from the equality in \re{fnkn}. $H(s,t)$ is as in \re{Hgen} and agrees with Sparre Andersen's formula \cite[Theorem 2]{andersenH_n}. 

\vspace{10pt}

\proof
Recall first that $H_n^{\rho} \ed H_n$ and $F_n^{\rho} \ed F_n$.
Let $n(s)$ be a geometric random variable with parameter $1-s$ and consider the walk of $n(s)$ steps. We have by definition
$$
H_{n(s)} = H_{s,0} + \sum_{j=1}^{\infty} H_{s,j}
$$
Thus the generating function of $H_{n(s)}$ is the product of the generating functions of $H_{s,0}$ and $H_{s,j}$, 
$j \ge 1$. These are known from Proposition \ref{prop:HKFdist}, thus
\begin{eqnarray*}
\sum_{m=0}^{\infty} t^m \P(H_{n(s)}=m) & = & e^{ (t-1) \mu_0(s) } \prod_{j=1}^{\infty} \frac{ e^{-\mu_j(s)} }{ 1-t + te^{-\mu_j(s)}} \\
& = & (1-s) e^{ t \mu_0(s) } \prod_{j=1}^{\infty} \frac{1}{ 1-t + te^{-\mu_j(s)}} 
\end{eqnarray*}
We can conclude that
\begin{eqnarray*}
H(s,t) & = & \sum_{n=0}^{\infty} \sum_{m=0}^{n} \P(H_n = m) s^n t^m \\
& = & (1-s)^{-1} \sum_{m=0}^{\infty} t^m \sum_{n=m}^{\infty} (1-s) s^n \P(H_n=m) \\
& = & (1-s)^{-1} \sum_{m=0}^{\infty} t^m \P(H_{n(s)}=m) \\
& = & e^{ t \mu_0(s) } \prod_{j=1}^{\infty} \frac{1}{ 1-t + te^{-\mu_j(s)}} 
\end{eqnarray*}
The deduction for $F(s,t)$ is similar, and as already mentioned, $K(s,t)$ is well known.
\endpf
\vspace{10pt}

In order to fully describe the joint law of the nested compositions, two more lemmas are necessary. The first contains information about the lengths of each segment or excursion, and the second describes how many excursions there are in each segment. We already know from the Poissonian description of the concave majorant the distribution of the number of segments with a given slope, and thus we already know the distribution of the number of segments within each face (see Theorem \ref{theorem:HKFcomposition} for the full description).

\begin{lemma}\label{lem:LHLF}
Consider the walk of $n(q)$ steps.
For $j \ge 1$, conditionally given $K_{q,j} = k_{q,j}$, let $L^K_{q,j,1},\ldots,L^K_{q,j,k_{q,j}}$ be the lengths of the $k_{q,j}$ segments of $S_{\rho}^{[0,n(q)]}$ of slope $x_j$. Then $L^K_{q,j,1}, \ldots , L^K_{q,j,k_{q,j}}$ are independent from each other and the lengths of all other segments> Moreover they are identically distributed with common probability generating function $G_{L_{q,j}^K}(z) = \mu_j(zq) / \mu_j(q)$.

For $j \ge 1$, conditionally given $H_{q,j} = h_{q,j}$, let $L^H_{q,j,1},\ldots,L^H_{q,j,h_{q,j}}$ be the lengths of the $h_{q,j}$ excursions of $S_{\rho}^{[0,n(q)]}$ of slope $x_j$. Then $L^K_{q,j,1}, \ldots , L^K_{q,j,h_{q,j}}$ are independent from each other and the lengths of all other segments. Moreover they are identically distributed with common probability generating function $G_{L_{q,j}^H}(z) = (1-e^{-\mu_j(zq)}) / (1-e^{-\mu_j(q)})$.

Furthermore, each excursion in the face of slope $x_j$ is independent and has the law of a random walk with increment distribution $X_1$ conditioned on making its first return to the line through the origin with slope $x_j$ before $n(q)$, an independent geometric random variable with parameter $1-q$, and remaining below that line before its first return time -- the excursion is taken to be that walk up to the time of its first return to the line with slope $x_j$.
\end{lemma}
\proof
By Poisson process properties, each $L^K_{q,j,1}, \ldots , L^K_{q,j,h_{q,j}}$ are independent from each other and the lengths of all other segments. By Poisson thinning, $\P(L^K_{q,j,1} = l) = l^{-1} q^l \P( S_k = kx_j)$, which gives the claimed generating function.

By the memoryless property of the geometric distribution of $n(q)$, each excursion of slope $x_j$ is independent, and is clearly independent from all excursions of other slopes. This gives the final assertion of the Lemma. By considering the total lengths of the face with slope $x_j$ we see that
$$
\sum_{i=1}^{H_{q,j}} L_{q,j,i}^H = \sum_{i=1}^{K_{q,j}} L_{q,j,i}^K
$$
By comparing the generating functions of both sides and using Proposition \ref{prop:HKFdist} we can deduce the claimed generating function $G_{L_{q,j}^H}(z)$.
\endpf
\vspace{10pt}

\begin{lemma}\label{lem:E_qj}
Conditionally given there are $k_{q,j}$ segments of $S_{\rho}^{[0,n(q)]}$ of slope $x_j$, let $E_{q,j,1},\ldots,E_{q,j,k_{q,j}}$ be the number of excursions in each of those $k_{n(q)}$ segments. Then $E_{q,j,1},\ldots,E_{q,j,k_{q,j}}$ are independent of each other and all other excursions and are identically distributed. Their common distribution is the log-series distribution with parameter $1-e^{-\mu_j(q)}$, that is
$$
\P( E_{q,j,1} = i) = \frac{(1-e^{-\mu_j(q)})^i}{i \mu_j(q)} \quad , i=1,2,\ldots 
$$
\end{lemma}
\proof
By Theorem \ref{theorem:poissonnoncont} the values of the walk $S_{\rho}^{[0,n(q)]}$ over each segment are independent, which gives the independence of $E_{q,j,1}, \ldots , E_{q,j,k_{q,j}}$.
By the independence of the excursions in the face of slope $x_j$ and the independence of the walks over each segment of slope $x_j$, $L^H_{q,j,1},\ldots,L^H_{q,j,E_{q,j}}$ are independent and identically distributed. By considering the total length of each segment of slope $x_j$, we have the identity in distribution
$$
L^K_{q,j,1} \ed \sum_{i=1}^{E_{q,j,1}} L^H_{q,j,1}
$$
which after applying generating function analysis reveals that
$$
G_{E_{q,j,1}}(z) : = \sum_{l=1}^{\infty} z^i \P ( E_{q,j,1} = i ) = \sum_{i=1}^{\infty} z^i \frac{(1-e^{-\mu_j(q)})^i}{i \mu_j(q)}
$$
\endpf
\vspace{10pt}

We are now ready to describe the joint law of the three nested compositions $\Xi_{[0,n(q)]}^H \preceq \Xi_{[0,n(q)]}^K \preceq \Xi_{[0,n(q)]}^F$. The following theorem is a summary of most of the information from Theorem \ref{theorem:poissonnoncont} to Lemma \ref{lem:E_qj}.
 
\begin{theorem}\label{theorem:HKFcomposition}
Let $n(q)$ be a geometric random variable with parameter $1-q$. Let $X_1,X_2,\ldots $ be independent and identically distributed. Let $S_j = \sum_{i=1}^j X_i$ for $j \ge 1$. Let $x_1,x_2,\ldots$ be an enumeration of the set of real numbers $x$ for which $\P(S_k=kx)$ is positive for some $k>0$, and for $j \ge 1$ let
$$
\mu_j(q) = \sum_{k=1}^{\infty} k^{-1} q^k \P(S_k=kx_j) 
$$
Let $S^{[0,n(q)]}_{\rho}$ be such that conditionally given $n(q)=n$, $S_{\rho}^{[0,n(q)]}$ is constructed in the same way as $S_{\rho}^{[0,n]}$ in Theorem \ref{theorem:1}. Let $\bar{C}_{\rho}^{[0,n(q)]}$ be the concave majorant of $S_{\rho}^{[0,n(q)]}$. Then independently for each $j \ge 1$:
\begin{itemize}
\item
There is a face of $\bar{C}_{\rho}^{[0,n(q)]}$  with slope $x_j$ with probability $1-e^{-\mu_j(q)}$.
\item
Conditionally given there is a face of slope $x_j$ the number of blocks of $\Xi_{[0,n]}^K$ with associated slope $x_j$ has the Poisson distribution with parameter $\mu_j(q)$, conditionally on the value being at least one.
\item
Conditionally given there are $k_{q,j}$ blocks of $\Xi_{[0,n]}^K$ with associated slope $x_j$, the number of excursions blocks in each of the $k_{q,j}$ segment blocks has the log-series distribution with parameter $1-e^{-\mu_j(q)}$, independently for each segment.
\item
The length of each excursion of slope $x_j$ is independent of all other excursions and has distribution with generating function 
$$
G_{L_{q,j}^H}(z) = (1-e^{-\mu_j(zq)}) / (1-e^{-\mu_j(q)})
$$
\end{itemize}
Any face block with associated slope $x$ such that $x \neq x_j$ for any $j \ge 1$ will be comprised of exactly one segment block, which will also be comprised of exactly one excursion block. The lengths and increments of faces with slope $x$ such that $x \neq x_j$ for any $j \ge 1$ form a Poisson point process on $\{ 1,2,\ldots \} \times \mathbb{R}$ with intensity $i^{-1} \P( S_i \in ds)$ for $i \ge 1, s \in \mathbb{R}$, but restricted to the region
$$
\{ (i,s) \in \{ 1,2,\ldots \} \times \mathbb{R} \, : \, s \neq i x_j \text{ for any } j \ge 1 \}
$$
Three nested compositions with the joint law of $\Xi_{[0,n(q)]}^H $, $\Xi_{[0,n(q)]}^K$ and $ \Xi_{[0,n(q)]}^F$ are created by uniformly randomly ordering the excursions within each segment, uniformly randomly ordering the segments within each face, arranging the faces in order of decreasing slope, and then looking at the induced compositions of excursion blocks, segment blocks and face blocks.
\end{theorem}

Theorem \ref{theorem:HKFcomposition} implies that  the compositions $\Xi_{[0,n(q)]}^H \preceq \Xi_{[0,n(q)]}^K \preceq \Xi_{[0,n(q)]}^F$ can be generated by nested renewal processes on $\mathbb{N}$ that terminate at some geometric time. There would be three types of renewal epochs. The first would be when a new face block started, which implies a new segment block and excursion block would also start. The second would be when only a new segment block and excursion block started, and the third would be when only a new excursion block started. Unlike in previous investigations into nested renewal sequences \cite{MR1733158,MR1395609}, the distributions of the length until the next renewal may change with time, and after a renewal has occurred, the number of future renewals may depend on how many have already occurred.

Theorem \ref{theorem:HKFcomposition} allows us to readily compute the probability of many fluctuation events for $S^{[0,n(q)]}$. Some examples are 

\begin{itemize}
\item
For each $j \ge 1$, the probability that $\bar{C}^{[0,n(q)]}$ consists of only one face of slope $x_j$ is $(1-q)^{-1} e^{-\mu_j(q)}$.
\item
The probability that $S^{[0,n(q)]}$ has a unique minimum, i.e.\ the probability that $\bar{C}^{[0,n(q)]}$ has no face of slope zero, is $\exp [ - \sum_{k=1}^{\infty} k^{-1} q^k \P(S_k=0) ]$.
\item
For each $j \ge 1$, the expected length of the face of $\bar{C}^{[0,n(q)]}$ of slope $x_j$ is $\sum_{k=1}^{\infty} q^k \P(S_k=kx_j)$.
\end{itemize}

\section{$S^{[0,n]}$ conditional on its concave majorant}
\label{sec:conditionedwalk}

To complete the rearrangement problem stated in the introduction, we now give a description of the law of $S^{[0,n]}$ conditional on  $\bar{C}^{[0,n]} = \bar{c}^{[0,n]}$. It is a generalization of the well known Vervaat transform for turning a bridge of a random walk into an excursion \cite[Theorem 5]{MR515820}.  
It relies on first choosing a segment composition $\Xi_{[0,n]}^K$ conditional on $\bar{C}_{\rho}^{[0,n]} = \bar{c}^{[0,n]}$ and then choosing a walk conditional on $\Xi_{[0,n]}^K$.

Let $\text{Supp}(\bar{C}^{[0,n]})$ be the support of the measure on concave functions on $[0,n]$ that represents the law of $\bar{C}^{[0,n]}$.
For any composition $(n_1,\ldots, n_k)$ of $n$ we say that $\sigma \in \Sigma_n$ is a $(n_1,\ldots,n_k)$-cyclic permutation of $[n]$ if its only action is to cyclically permute the first $n_1$ elements of $[n]$, cyclically permute the next $n_2$ elements of $[n]$ and so on. For example, $234175689$ is a $(4,3,2)$-cyclic permutation of $[9]$. Recall that in Section \ref{sec:noncontinuous} we defined $\mathcal{N}_n$ to be the set of compositions of $n$, and $\mathcal{N}(\bar{c}^{[0,n]}) \subseteq \mathcal{N}_n$ to be the set of possible values of $\Xi_{[0,n]}^K$ conditionally given $\bar{C}_{\rho}^{[0,n]} = \bar{c}^{[0,n]}$.

\begin{theorem}\label{theorem:conditionedwalk}
Let $S_0 = 0$ and $S_j = \sum_{\ell =1}^j X_{\ell}$ for $1 \le j \le n$, where $X_1,\ldots,X_n$ 
are exchangeable random variables. Let $S^{[0,n]} = \{ (j,S_j) \: : \: 0 \le j \le n \}$ and let $\bar{C}^{[0,n]}$ be the concave majorant of $S^{[0,n]}$. 
Suppose $\bar{c}^{[0,n]} \in \text{Supp}(\bar{C}^{[0,n]})$.
Let $q(\cdot)$ be the probability density function on $\mathcal{N}_n$ that is the regular conditional distribution of $\Xi_{[0,n]}$ conditionally given $\bar{C}_{\rho}^{[0,n]} = \bar{c}^{[0,n]}$.  
Let  $(N_{n,1}, N_{n,2} , \ldots ,  N_{n,K_n})$ be a composition of $n$ chosen according to the density function $q(\cdot)$, independently of $\{ X_j : 1 \le j \le n \}$. 

Conditionally given $\{ K_n = k \}$ and $\{ N_{n,i} = n_i \: : \: 1 \le i \le k \}$, let $Y_1, \ldots , Y_n$ be random variables, independent of all previously introduced random variables, whose joint law that is the regular conditional joint distribution of $X_1, \ldots , X_n$ conditionally given $\{  S_j \in d\bar{c}^{[0,n]}(j) , j=\sum_{i=1}^m n_i, 1 \le m \le k \}$.

Conditionally given $Y_1,\ldots,Y_n$, let $B$ be the random set of $(n_1,\ldots,n_k)$-cyclic permutations of $[n]$ such that
$$
Y_{\sigma(j)} \geq \bar{c}^{[0,n]}(j) \quad \text{ for } 1 \le j \le n
$$
if and only if $\sigma \in B$. Let $\hat{\rho}$ be an independently chosen uniform random element of $B$, and let $S^{\hat{\rho}}_j = \sum_{\ell=1}^j Y_{\hat{\rho} (\ell)}$ for $1 \le j \le n$. Then $S_{\hat{\rho}}^{[0,n]} := \{ (j,S_j^{ \hat{\rho}} ) \: : \: 1 \le j \le n \}$ has the regular conditional distribution of $S^{[0,n]}$ conditionally given $\bar{C}^{[0,n]} = \bar{c}^{[0,n]}$.
\end{theorem}

The theorem is direct result of Bayes' rule and Theorem \ref{theorem:1}. Note that when $X_1,\ldots,X_n$ satisfy assumption $\mathbf{A}$, $\mathcal{N}(\bar{c}^{[0,n]}) $ has only one element, the composition induced by the lengths of the faces of $\bar{c}^{[0,n]}$, and $A$ also only contains one element by Lemma \ref{lem:cyclics}, so the theorem simplifies significantly. It remains to describe $q(\cdot)$.

\begin{lemma}\label{lem:conditionalcomposition}
Suppose $\bar{c}^{[0,n]} \in \text{Supp}(\bar{C}^{[0,n]})$ and that $X_1,\ldots,X_n$ are exchangeable. The regular conditional distribution of $\Xi_{[0,n]}$ conditionally given $\bar{C}_{\rho}^{[0,n]} = \bar{c}^{[0,n]}$ is given by
$$
 \P( \bar{C}^{[0,n]}(j) \in d\bar{c}^{[0,n]}(j) , 1 \le j \le n ) \P( \Xi^K_{[0,n]} = (n_1,\ldots,n_k) | \bar{C}_{\rho}^{[0,n]} = \bar{c}^{[0,n]} ) 
$$
$$
 = 1_{ (n_1,\ldots,n_k) \in \mathcal{N}(\bar{c}^{[0,n]}) }
\frac{\prod_{i=1}^k n_i }{\prod_{j=1}^{f(\bar{c}^{[0,n]} )} k_j(n_1,\ldots,n_k)!} 
\P ( S_j \in d\bar{c}^{[0,n]}(j) , j= \mbox{$\sum_{i=1}^l n_i$}, 1 \le l \le k )
$$
where $S_j$, $ 1 \le j \le n$ is as in Theorem \ref{theorem:conditionedwalk}.
\end{lemma}
\proof
Let $(n_1,\ldots,n_k) \in \mathcal{N}(\bar{c}^{[0,n]}) $. Following the construction in Theorem \ref{theorem:1}, by the Ewens sampling formula  the probability that $\{ L_{n,1}, \ldots, L_{n,K_n} \}$ is a list of the elements of $(n_1,\ldots,n_k)$ in non-increasing order is $\left( \prod_{j=1}^n (a_j!)^{-1} \right) \left( \prod_{i=1}^k n_i^{-1} \right)$ where $a_j = \# \{ i \: : \: 1 \le i \le k , n_i = j \}$ for $1 \le j \le n$. Conditionally given  $\{ L_{n,1}, \ldots, L_{n,K_n} \}$ is a list of the elements of $(n_1,\ldots,n_k)$ in non-increasing order the probability of the event 
$ \{ \Xi^K = (n_1,\ldots,n_k) , \bar{C}^{[0,n]} = \bar{c}^{[0,n]} \}$ is 
$$
\left( \frac{\prod_{j=1}^n a_j! }{\prod_{j=1}^{f(\bar{c}^{[0,n]} )} k_j(n_1,\ldots,n_k)!} \right)
\P ( S_j \in d\bar{c}^{[0,n]}(j) , j=\sum_{i=1}^l n_i, 1 \le l \le k )
$$
where the denominator in the multiplicative factor in the brackets is due to the restrictions on the orderings of partition segments within each face, and the numerator is because of repeated segment lengths.
\endpf
\vspace{10pt}

We say that the concave majorant of a walk is \emph{trivial} if it has only one face. A particularly useful form of Theorem \ref{theorem:conditionedwalk} arises from the special case when the increments $X_1,\ldots,X_n$ are independent, the probability that the concave majorant of $S^{[0,n]}$ is trivial with slope zero is positive, and we want the conditional distribution of the walk $S^{[0,n]}$ given it has trivial concave majorant of slope zero. 
By subtraction of a line of constant slope, this gives us the conditional distribution of the walk $S^{[0,n]}$ given it has trivial concave majorant of any slope, as long as the probability that the concave majorant of $S^{[0,n]}$ is trivial with that slope is positive. In the case where we want the regular conditional distribution for $S^{[0,n]}$ conditional on having trivial concave majorant of a slope that has zero probability, then the only possible value for $\Xi_{[0,n]}$ is the trivial composition $(n)$.

\begin{corollary}\label{cor:flatwalk}
Let $S_0 = 0$ and $S_j = \sum_{i=1}^j X_i$ for $1 \le j \le n$, where $X_1, \ldots,X_n$ 
are independent identically distributed random variables, and let $S^{[0,n]} = \{ (j,S_j) \: : \: 0 \le j \le n \}$.
Suppose that 
$$
p_{\text{triv}} := \P( \text{concave majorant of $S^{[0,n]}$ is trivial with slope zero} ) > 0
$$
Define a probability density function $q(\cdot)$ on $\mathcal{N}_n$ by 
$$
q \left( (n_1,\ldots,n_k) \right) = 
\frac{1}{p_{\text{triv}} k!}  
\prod_{i=1}^k n_i u_{n_i} 
$$
where $u_j = \P(S_j=0)$ for $1 \le j \le n$. Let  $(N_{n,1}, N_{n,2} , \ldots ,  N_{n,K_n})$ be a composition of $n$ chosen according to the density function $q(\cdot)$, independently of $\{ X_j : 1 \le j \le n \}$. 

Conditionally given $\{ K_n = k \}$ and $\{ N_{n,i} = n_i \: : \: 1 \le i \le k \}$, independently for each $1 \le i \le k$ let $Y_{n_1+\cdots+n_{i-1}+1},\ldots,Y_{n_1+\cdots+n_i}$ be random variables, independent of all previously introduced random variables, whose joint law that is the regular conditional joint distribution of $X_1, \ldots , X_{n_i}$ conditionally given $\sum_{\ell=1}^{n_i} X_{\ell} = 0$.

Conditionally given $Y_1,\ldots,Y_n$, let $B$ be the random set of $(n_1,\ldots,n_k)$-cyclic permutations of $[n]$ such that
$$
Y_{\sigma(j)} \leq 0 \quad \text{ for } 1 \le j \le n
$$
if and only if $\sigma \in B$. Let $\hat{\rho}$ be an independently chosen uniform random element of $B$, and let $S^{\hat{\rho}}_j = \sum_{\ell=1}^j Y_{\hat{\rho} (\ell)}$ for $1 \le j \le n$. Then $S_{\hat{\rho}}^{[0,n]} := \{ (j,S_j^{ \hat{\rho}} ) \: : \: 1 \le j \le n \}$ has the regular conditional distribution of $S^{[0,n]}$ conditionally given $S^{[0,n]}$ has trivial concave majorant with slope zero.
\end{corollary}

\section{A path transformation}
\label{sec:pathtransform}

This section provides an important path transformation which by taking scaling limits is used by Pitman and Uribe Bravo to completely describe the concave majorant (or as in that paper, convex minorant) of a L\'{e}vy process and the excursions of that process beneath its concave majorant \cite{pitmanbravo}. Essentially, the idea is that a uniformly sampled face of the concave majorant should have uniform length and the walk over it should be a Vervaat like transform of some walk of the same length.

Let $S_0 = 0$ and $S_j = \sum_{i=1}^n X_i$ for $1 \le j \le n$, where $X_i, i = 1, \ldots , n$ 
are exchangeable  random variables satisfying assumption $\mathbf{A}$.
We introduce the following path transformation
for the random walk $S^{[0,n]} = \{ (j,S_j) \: , \: 1 \le j \le n \}$. Let $U$ be distributed uniformly on $[n]$. Let $g$ and $d$ be the left and right end points respectively of the face of the concave majorant of $S^{[0,n]}$ containing the $U$th increment $X_U$. Define $S_j^U $ for $1 \le j \le n $ by
\eq
\lb{pathtrans}
S^U_j =\begin{cases}
S_{U+j}-S_{U}& \text{for $0 \le j < d-U$} \\
S_{g+j-(d-U)}+S_d-S_g-S_U&\text{for $d - U \le j < d-g$}\\
S_{j-(d-g)} +S_d-S_g&\text{for $d-g \le j < d$}\\
S_j&\text{for $d \le j \le n$}.
\end{cases}
\en 
and let $S^{[0,n]}_U = \{ (j,S_j^U) \: , \: 1 \le j \le n \}$.

\begin{figure}
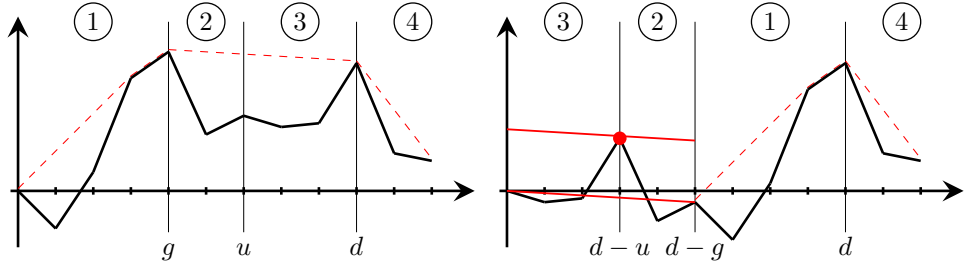

\begin{pgfpicture}{0cm}{0cm}{9cm}{5.02cm}
\begin{pgfscope} 
\pgfsetendarrow{\pgfarrowsingle}
\pgfsetlinewidth{1pt} 
\pgfxyline(-0.1,0.75)(6,0.75) 
\pgfxyline(0,0)(0,3) 
\pgfxyline(6.4,0.75)(12.5,0.75) 
\pgfxyline(6.5,0)(6.5,3) 
\end{pgfscope} 
\begin{pgfscope} 
\pgfsetlinewidth{1pt} 
\pgfxyline(0.5,0.7)(0.5,0.8) 
\pgfxyline(1.0,0.7)(1.0,0.8) 
\pgfxyline(1.5,0.7)(1.5,0.8) 
\pgfxyline(2.0,0.7)(2.0,0.8) 
\pgfxyline(2.5,0.7)(2.5,0.8) 
\pgfxyline(3.0,0.7)(3.0,0.8) 
\pgfxyline(3.5,0.7)(3.5,0.8) 
\pgfxyline(4.0,0.7)(4.0,0.8) 
\pgfxyline(4.5,0.7)(4.5,0.8) 
\pgfxyline(5.0,0.7)(5.0,0.8) 
\pgfxyline(5.5,0.7)(5.5,0.8) 
\pgfxyline(7.0,0.7)(7.0,0.8) 
\pgfxyline(7.5,0.7)(7.5,0.8) 
\pgfxyline(8.0,0.7)(8.0,0.8) 
\pgfxyline(8.5,0.7)(8.5,0.8) 
\pgfxyline(9.0,0.7)(9.0,0.8) 
\pgfxyline(9.5,0.7)(9.5,0.8) 
\pgfxyline(10.0,0.7)(10.0,0.8) 
\pgfxyline(10.5,0.7)(10.5,0.8) 
\pgfxyline(11.0,0.7)(11.0,0.8) 
\pgfxyline(11.5,0.7)(11.5,0.8) 
\pgfxyline(12.0,0.7)(12.0,0.8) 
\end{pgfscope} 
\begin{pgfscope}
\pgfsetlinewidth{0.1pt}
\pgfxyline(2,0.2)(2,3) 
\pgfxyline(3,0.2)(3,3) 
\pgfxyline(4.5,0.2)(4.5,3) 
\pgfputat{\pgfxy(2,-0.1)}{\pgfbox[center,base]{$g$}}
\pgfputat{\pgfxy(3,-0.1)}{\pgfbox[center,base]{$u$}}
\pgfputat{\pgfxy(4.5,-0.1)}{\pgfbox[center,base]{$d$}}
\pgfxyline(8,0.2)(8,3) 
\pgfxyline(9,0.2)(9,3) 
\pgfxyline(11,0.2)(11,3) 
\pgfputat{\pgfxy(8,-0.1)}{\pgfbox[center,base]{$d-u$}}
\pgfputat{\pgfxy(9,-0.1)}{\pgfbox[center,base]{$d-g$}}
\pgfputat{\pgfxy(11,-0.1)}{\pgfbox[center,base]{$d$}}
\end{pgfscope}
\begin{pgfscope}
\color{red}
\pgfsetdash{{3pt}{3pt}}{0pt}
\pgfxyline(0.0,0.78)(1.5,2.28)
\pgfxyline(1.5,2.28)(2.0,2.63) 
\pgfxyline(2.0,2.63)(4.5,2.48) 
\pgfxyline(4.5,2.48)(5.5,1.18)
\pgfxyline(9.0,0.63)(10.5,2.13) 
\pgfxyline(10.5,2.13)(11.0,2.48) 
\pgfxyline(11.0,2.48)(12.0,1.18)
\color{black}
\end{pgfscope}
\begin{pgfscope} 
\pgfsetlinewidth{1pt} 
\pgfxyline(0,0.75)(0.5,0.25) 
\pgfxyline(0.5,0.25)(1.0,1.0) 
\pgfxyline(1.0,1.0)(1.5,2.25) 
\pgfxyline(1.5,2.25)(2,2.6) 
\pgfxyline(2,2.6)(2.5,1.5) 
\pgfxyline(2.5,1.5)(3,1.75) 
\pgfxyline(3,1.75)(3.5,1.6) 
\pgfxyline(3.5,1.6)(4,1.65) 
\pgfxyline(4,1.65)(4.5,2.45) 
\pgfxyline(4.5,2.45)(5,1.25) 
\pgfxyline(5,1.25)(5.5,1.15) 
\pgfxyline(6.5,0.75)(7,0.6) 
\pgfxyline(7.0,0.6)(7.5,0.65) 
\pgfxyline(7.5,0.65)(8,1.45) 
\pgfxyline(8,1.45)(8.5,0.35) 
\pgfxyline(8.5,0.35)(9,0.6) 
\pgfxyline(9,0.6)(9.5,0.1) 
\pgfxyline(9.5,0.1)(10,0.85) 
\pgfxyline(10,0.85)(10.5,2.1) 
\pgfxyline(10.5,2.1)(11,2.45) 
\pgfxyline(11,2.45)(11.5,1.25) 
\pgfxyline(11.5,1.25)(12,1.15) 
\pgfputat{\pgfxy(1,3)}{\pgfbox[center,center]{$1$}}
\pgfputat{\pgfxy(2.5,3)}{\pgfbox[center,center]{$2$}}
\pgfputat{\pgfxy(3.75,3)}{\pgfbox[center,center]{$3$}}
\pgfputat{\pgfxy(5.25,3)}{\pgfbox[center,center]{$4$}}
\pgfsetlinewidth{0.5pt}
\pgfcircle[stroke]{\pgfxy(5.25,3)}{0.25cm}
\pgfcircle[stroke]{\pgfxy(3.75,3)}{0.25cm}
\pgfcircle[stroke]{\pgfxy(2.5,3)}{0.25cm}
\pgfcircle[stroke]{\pgfxy(1,3)}{0.25cm}
\pgfsetlinewidth{1pt}
\pgfputat{\pgfxy(7.25,3)}{\pgfbox[center,center]{$3$}}
\pgfputat{\pgfxy(8.5,3)}{\pgfbox[center,center]{$2$}}
\pgfputat{\pgfxy(10,3)}{\pgfbox[center,center]{$1$}}
\pgfputat{\pgfxy(11.75,3)}{\pgfbox[center,center]{$4$}}
\pgfsetlinewidth{0.5pt}
\pgfcircle[stroke]{\pgfxy(7.25,3)}{0.25cm}
\pgfcircle[stroke]{\pgfxy(8.5,3)}{0.25cm}
\pgfcircle[stroke]{\pgfxy(10,3)}{0.25cm}
\pgfcircle[stroke]{\pgfxy(11.75,3)}{0.25cm}
\color{red}
\pgfcircle[fill]{\pgfxy(8,1.45)}{0.09cm}
\pgfsetlinewidth{0.7pt}
\color{red}
\pgfxyline(6.5,0.75)(9.0,0.60)
\color{red}
\pgfxyline(6.5,1.57)(9.0,1.42)
\end{pgfscope} 
\end{pgfpicture}
\captionsetup{singlelinecheck=off}
\caption{An example of the ``$3214$'' path transformation of Theorem
  \ref{pathtransform}. The walk on the right is the transformed version
of the walk on the left. Note how given $d-g$ the transform is easily
inverted - the index at which the first $d-g$ increments should start
after cyclic permutation is marked, and can be found by lowering a line
with the slope the mean of the first $d-g$ increments.}
\label{fig:transform}
\end{figure}

\begin{theorem}
\label{pathtransform}
$$
(U,S^{[0,n]}) \ed (d-g,S^{[0,n]}_U)
$$
\end{theorem}

\vspace{10pt}

In fact, Theorem \ref{pathtransform} provides an alternative method of proving Theorem \ref{theorem:1} under assumption $\mathbf{A}$, since by applying the transformation again to the $S_U^{[0,n]}$ restricted to the interval $[d-g,n]$, and then doing this repeatedly until there is nothing left to transform, we are actually performing the inverse of the transformation given in Theorem \ref{theorem:1}. However, this method does not extend to cover the general case as considered in Section \ref{sec:noncontinuous}, so we will not expand on it.

\vspace{10pt}

\proof
As in the proof of Theorem \ref{theorem:1} under assumption $\mathbf{A}$ in Section \ref{sec:proof}, it is enough to show that the equality in distribution holds when $X_1,\ldots,X_n$ are samples without replacement from $x_1,\ldots,x_n$ satisfying assumption $\mathbf{A}$. $S^{[0,n]}$ and $S^{[0,n]}_U$ may thus be thought of as permutations of $n$, so we may think of the mapping $(U,S^{[0,n]}) \mapsto (d-g,S^{[0,n]}_U)$ as a mapping from $[n] \times \Sigma_n$ to itself. Since $U$ is uniform on $[n]$, and the ordering of $X_1,\ldots,X_n$ is a uniform random permutation of $x_1,\ldots,x_n$, it is enough to show that this mapping is a bijection. To do this, it suffices to show that the mapping is surjective. This can be seen visually in Figure \ref{fig:transform} since it is clear from the figure and its description that the map is easily inverted. More formally, to show that the map is surjective it is sufficient to show that for $k \in [n]$ there exists $u \in [n]$ and $\sigma \in \Sigma_n$ such that 
$$
 \left( u, \{ (0,0),(1,x_{\sigma(1)}),(2,x_{\sigma(1)} +x_{\sigma_2}), \ldots, (n, \sum_{i=1}^n x_{\sigma(i)}) \} \right) \: \: \: \:  \: \: \: \:  \: \: \: \: \: \: \: \:  \: \: \: \:  \: \: \: \:  \: \: \: \:
$$
$$
 \: \: \: \:  \: \: \: \:  \: \: \: \: \: \: \: \:  \: \: \: \:  \: \: \: \: \mapsto 
\left( k , \{ (0,0),(1,x_1),(2,x_1 +x_2), \ldots, (n, \sum_{i=1}^n x_i) \} \right)
$$
Let $f$ be the number of faces of the concave majorant of the walk of length $n-k$ with increments $x_{k+1},\ldots,x_n$, and let the lengths and increments of these faces in order of appearance be $(\ell_1,s_1),\ldots,(\ell_f,s_f)$. Let $r$ be the unique $r \in [k]$ such that the walk with increments 
$$(x_{r+1},x_{(r+1)\! \! \! \! \! \mod k \: + 1 },x_{(r+2)\! \! \! \! \! \mod k \: + 1 }, \ldots,x_{(r+k-2)\! \! \! \! \! \mod k \: + 1 } , x_r)$$
remains below its concave majorant. Let $s^* = \sum_{i=1}^k x_i$, and let $m$ be the unique $m \in \{ 0,\ldots,f \}$ such that
$$  \frac{s_m}{\ell_m} > \frac{s^*}{k} > \frac{s_{m+1}}{\ell_{m+1}}  $$
where we say that $s_0/\ell_0 = + \infty$ and $s_{f+1}/\ell_{f+1} = \infty$. The appropriate $(\sigma(1),\ldots,\sigma(n))$ is given by
$$
\begin{array}{l}
( k+1,k+2,\ldots, k+ \mbox{$\sum_{i=1}^m \ell_i$} , \\
\\
 \: \: \:  \: \: \: \:  \: \: 
r+1, (r+1)\! \! \! \! \! \mod k \: + 1 , (r+2)\! \! \! \! \! \mod k \: + 1, \ldots , r, 
 k+\mbox{$\sum_{i=1}^m \ell_i$}  + 1, \ldots, n)
\end{array}
$$
\endpf
\vspace{10pt}

\bibliographystyle{plain} 
\bibliography{refs.bib} {}

\begin{thebibliography}{10}

\bibitem{andersenH_n}
Erik~Sparre Andersen.
\newblock On the distribution of the random variable {$H_n$}. {T}ech. {S}ci.
  {N}ote {N}o. 1, {C}ontract {N}o. {AF} 61(052)-42, {F}ebruary 27, 1959.

\bibitem{MR0058893b}
Erik~Sparre Andersen.
\newblock On the fluctuations of sums of random variables {II}.
\newblock {\em Math. Scand.}, 2:195--223, 1954.

\bibitem{MR1733158}
Jean Bertoin.
\newblock Renewal theory for embedded regenerative sets.
\newblock {\em Ann. Probab.}, 27(3):1523--1535, 1999.

\bibitem{MR1747095}
Jean Bertoin.
\newblock The convex minorant of the {C}auchy process.
\newblock {\em Electron. Comm. Probab.}, 5:51--55 (electronic), 2000.

\bibitem{MR0162302}
H.~D. Brunk.
\newblock A generalization of {S}pitzer's combinatorial lemma.
\newblock {\em Z. Wahrscheinlichkeitstheorie und Verw. Gebiete}, 2:395--405
  (1964), 1964.

\bibitem{MR1395609}
Persi Diaconis, Susan Holmes, Svante Janson, Steven~P. Lalley, and Robin
  Pemantle.
\newblock Metrics on compositions and coincidences among renewal sequences.
\newblock In {\em Random discrete structures ({M}inneapolis, {MN}, 1993)},
  volume~76 of {\em IMA Vol. Math. Appl.}, pages 81--101. Springer, New York,
  1996.

\bibitem{MR0270403}
William Feller.
\newblock {\em An introduction to probability theory and its applications.
  {V}ol. {II}.}
\newblock Second edition. John Wiley \& Sons Inc., New York, 1971.

\bibitem{gp03}
A.~Gnedin and J.~Pitman.
\newblock Regenerative composition structures.
\newblock {\em Ann. Probab.}, 33(2):445--479, 2005.

\bibitem{MR994088}
Charles~M. Goldie.
\newblock Records, permutations and greatest convex minorants.
\newblock {\em Math. Proc. Cambridge Philos. Soc.}, 106(1):169--177, 1989.

\bibitem{MR588409}
Priscilla Greenwood and Jim Pitman.
\newblock Fluctuation identities for {L}\'evy processes and splitting at the
  maximum.
\newblock {\em Adv. in Appl. Probab.}, 12(4):893--902, 1980.

\bibitem{greenpit80}
Priscilla Greenwood and Jim Pitman.
\newblock Fluctuation identities for random walk by path decomposition at the
  maximum.
\newblock {\em Advances in Applied Probability}, 12(2):291--293, 1980.

\bibitem{MR714964}
Piet Groeneboom.
\newblock The concave majorant of {B}rownian motion.
\newblock {\em Ann. Probab.}, 11(4):1016--1027, 1983.

\bibitem{MR0062867}
M.~Kac.
\newblock Toeplitz matrices, translation kernels and a related problem in
  probability theory.
\newblock {\em Duke Math. J.}, 21:501--509, 1954.

\bibitem{MR2245368}
J.~Pitman.
\newblock {\em Combinatorial stochastic processes}, volume 1875 of {\em Lecture
  Notes in Mathematics}.
\newblock Springer-Verlag, Berlin, 2006.
\newblock Lectures from the 32nd Summer School on Probability Theory held in
  Saint-Flour, July 7--24, 2002, With a foreword by Jean Picard.

\bibitem{pitmanbravo}
J.~{Pitman} and G.~{Uribe Bravo}.
\newblock {The convex minorant of a L\'evy process}.
\newblock {\em Ann. Probab.}, 2011.
\newblock To appear.

\bibitem{MR2205729}
Zhihua Qiao and J.~Michael Steele.
\newblock Random walks whose concave majorants often have few faces.
\newblock {\em Statist. Probab. Lett.}, 75(2):97--102, 2005.

\bibitem{MR0195117}
L.~A. Shepp and S.~P. Lloyd.
\newblock Ordered cycle lengths in a random permutation.
\newblock {\em Trans. Amer. Math. Soc.}, 121:340--357, 1966.

\bibitem{MR0169306}
S.~Sherman.
\newblock Fluctuation and periodicity.
\newblock {\em J. Math. Anal. Appl.}, 9:468--476, 1964.

\bibitem{MR0079851}
Frank Spitzer.
\newblock A combinatorial lemma and its application to probability theory.
\newblock {\em Trans. Amer. Math. Soc.}, 82:323--339, 1956.

\bibitem{MR1910531}
J.~Michael Steele.
\newblock The {B}ohnenblust-{S}pitzer algorithm and its applications.
\newblock {\em J. Comput. Appl. Math.}, 142(1):235--249, 2002.
\newblock Probabilistic methods in combinatorics and combinatorial
  optimization.

\bibitem{MR515820}
Wim Vervaat.
\newblock A relation between {B}rownian bridge and {B}rownian excursion.
\newblock {\em Ann. Probab.}, 7(1):143--149, 1979.

\end{thebibliography}
\end{document}